\documentclass[10pt, a4paper, reqno]{amsart}
\usepackage{amssymb}
\usepackage[arrow, matrix, curve]{xy}
\usepackage{esint}
\usepackage[shortlabels, inline]{enumitem}
\usepackage{url}
\usepackage{tikz}
\usepackage{tikz-cd}
\usepackage{stmaryrd}
\usepackage{geometry}
\usepackage[numeric, abbrev, nobysame]{amsrefs}
\usepackage{color}

\makeatletter
\@namedef{subjclassname@2020}{\textup{2020} Mathematics Subject Classification}
\makeatother

\theoremstyle{plain}
\newtheorem{thm}{Theorem}
\newtheorem{prop}[thm]{Proposition}
\newtheorem{lemma}[thm]{Lemma}
\newtheorem{cor}[thm]{Corollary}

\let\frak\mathfrak

\theoremstyle{definition}
\newtheorem{defin}{Definition}

\theoremstyle{remark}
\newtheorem*{rem}{Remark}

\DeclareMathOperator{\SO}{SO}

\DeclareMathOperator{\GL}{GL}
\DeclareMathOperator{\Hom}{Hom}
\DeclareMathOperator{\re}{re}
\DeclareMathOperator{\im}{im}

\DeclareMathOperator{\vol}{vol}

\DeclareMathOperator{\ad}{ad}

\DeclareMathOperator{\SU}{SU}
\DeclareMathOperator{\U}{U}

\DeclareMathOperator{\CSU}{CSU}

\newcommand{\x}{\times}

\newcommand{\mc}{\mathbb{C}}
\newcommand{\mE}{\mathbb{E}}

\newcommand{\mh}{\mathbb{H}}
\newcommand{\mk}{\mathbb{K}}

\newcommand{\mr}{\mathbb{R}}

\newcommand{\mv}{\mathbb{V}}
\newcommand{\mw}{\mathbb{W}}
\newcommand{\mz}{\mathbb{Z}}

\newcommand{\al}{\alpha}
\newcommand{\be}{\beta}

\renewcommand{\o}{\circ}

\newcommand{\xa}{\mathfrak{a}}
\newcommand{\xg}{\mathfrak{g}}

\newcommand{\xk}{\mathfrak{k}}
\newcommand{\xm}{\mathfrak{m}}
\newcommand{\xn}{\mathfrak{n}}
\newcommand{\xp}{\mathfrak{p}}
\newcommand{\xq}{\mathfrak{q}}

\newcommand{\cc}{\mathcal{C}}
\newcommand{\cd}{\mathcal{D}}
\newcommand{\ocd}{\overline{\mathcal{D}}}
\newcommand{\ce}{\mathcal{E}}

\newcommand{\ch}{\mathcal{H}}
\newcommand{\cI}{\mathcal{I}}
\newcommand{\cK}{\mathcal{K}}

\newcommand{\cp}{\mathcal{P}}

\newcommand{\cl}{\mathcal{L}}

\newcommand{\cU}{\mathcal{U}}

\newcommand{\bz}{\overline{z}}

\newcommand{\opartial}{\overline{\partial}}
\newcommand{\oBox}{\overline{\Box}}

\newcommand{\otau}{\overline{\tau}}
\newcommand{\uPhi}{\underline{\Phi}}

\begin{document}

\title{Poisson transforms, the BGG complex, and \\
  discrete series representations of $SU(n+1,1)$.}

\author{Andreas \v{C}ap, Christoph Harrach, and Pierre Julg}

\address{A.\ \v C.\ and C.\ H.: Faculty of Mathematics\\
  University of Vienna\\
  Oskar-Morgenstern Platz 1\\
  1090 Wien\\
  Austria\\
P.\ J.:  MAPMO UMR 7349, F\'{e}d\'{e}ration Denis Poisson \\
Universit\'{e} d'Orl\'{e}ans, Coll\'{e}gium Sciences et Techniques \\
B\^{a}timent de math\'{e}\-ma\-tiques - Route de Chartres \\
B.P. 6759, 45067 Orl\'{e}ans cedex 2 \\
France}
\email{Andreas.Cap@univie.ac.at, Christoph.Harrach@univie.ac.at, pierre.julg@univ-orleans.fr}


\begin{abstract}
The aim of this article is to construct specific Poisson transforms mapping
differential forms on the sphere $S^{2n+1}$ endowed with its natural CR structure to
forms on complex hyperbolic space. These transforms have co-closed harmonic values,
descend to the BGG (Rumin) complex, and intertwine the differential operators in that 
complex with the exterior derivative.

Passing to the Poincar\'e ball model, we analyze boundary asymptotics, proving that 
the values of our transforms admit a continuous extension to the boundary in degrees 
$\leq n$. Finally, we show that composing the exterior derivative with the transform 
in degree $n$, one obtains an isomorphism between the kernel of the Rumin operator in 
degree $n$ and a dense subspace of the $L^2$-harmonic forms on complex hyperbolic space. 
This provides a realization of the discrete series representations of $SU(n+1,1)$ with 
trivial infinitesimal character in spaces of differential forms on the compact manifold 
$S^{2n+1}$.

These developments are motivated by a program of the third author to prove some instances 
of the Baum-Connes conjecture. The first part of the article is valid in a much more 
general setting, in particular, it is relevant for cases in which the conjecture
is still open. 
\end{abstract}

\maketitle

\pagestyle{myheadings} \markboth{\v Cap, Harrach, Julg}{Poisson transforms and
  discrete series representations}

\section{Introduction}
The idea to extend the Poisson transform to an operation defined on differential
forms goes back to the work \cite{gaillard} of P.Y.\ Gaillard. He defined a transform
on differential forms on a sphere with values in real hyperbolic space of one higher
dimension in order to construct differential forms on hyperbolic space which are
coclosed and harmonic. This transform is equivariant for the natural actions of
$SO_0(n+1,1)$ by conformal transformations on $S^n$ and by isometries of the
hyperbolic metric on hyperbolic space and representation theory plays an important
role in the construction. 

  It is natural to study the asymptotics of the values of the transform towards
  infinity and it turns out this depends crucially on the degree of differential
  forms. For sufficiently small degree, these values admit a smooth extension to the
  boundary, while for high degrees they diverge towards the
    boundary. Particularly interesting behavior occurs in the case that $n$ is odd
  for forms whose degree is exactly half the dimension of the interior hyperbolic
  space. In this case, the transform leads to $L^2$-harmonic forms on hyperbolic
  space, which only exist in this specific degree. The subspace of $L^2$-harmonic
  forms is invariant under the action of $SO_0(n+1,1)$ and realizes the direct sum of
  the two discrete series representations of this group with trivial infinitesimal
  character. Building on Gaillard's work, J.\ Lott showed in \cite{lott} that these
  discrete series representations can be identified with spaces of closed forms on
  the boundary sphere of Sobolev regularity $H^{-1/2}$.

  There were several attempts to generalize these methods to the other symmetric
  spaces of real rank one, i.e.\ the complex and quaternionic hyperbolic spaces and
  the octonionic hyperbolic plane. But it quickly turned out that the real case is
  deceivingly simple. For $SO_0(n+1,1)$ there is essentially only one possible
  transform in each degree, which makes it easy to relate the transform to natural
  differential operators in the interior. Already in the complex case, corresponding
  to $\SU(n+1,1)$, there are many different choices for Poisson transforms on
  differential forms, not all of them produce harmonic values, and so on. The
  question of harmonic values turns out to be very nicely related to the geometry of
  the sphere in each case. In the complex case, which is mainly considered in this
  article, one obtains the sphere $S^{2n+1}$ with the CR structure coming from the
  embedding into $\mc^{n+1}$. This admits a natural interpretation as the structure
  on the boundary at infinity of complex hyperbolic space, which is a complete
  K\"ahler-Einstein manifold. In the quaternionic case, one obtains a sphere
  $S^{4n+3}$ endowed with a so-called quaternionic contact structure which comes from
  the embedding into $\mathbb H^{n+1}$. This structure can be interpreted as coming
  from being a boundary at infinity of quaternionic hyperbolic space which is a
  complete quaternion-K\"ahler manifold, see \cite{biquard} for both the complex and
  the quaternionic case.

  The geometries on the spheres one obtains all fall into the class of so-called
  parabolic geometries, which makes a large set of efficient tools for their study
  available, see \cite{cap_slovak}. Except for the real case, this structure leads to
  a refinement of the de Rham complex, which in the complex case is known as the
  Rumin complex. Writing the sphere as $G/P$ for the appropriate
    semi-simple group $G$ and parabolic subgroup $P\subset G$, the representations
  of $P$ which induce the bundles of differential forms are not completely reducible
  in this case, but one can form completely reducible subquotients in a natural
  way. These induce subquotient bundles and there is a construction of natural
  differential operators (partly of higher order) acting between sections of
  those. This can be done in such a way that one obtains a complex that computes the
  de Rham cohomology. This is a very special case of the construction of
  Bernstein-Gelfand-Gelfand (BGG) sequences, whence we refer to the resulting
  complexes as BGG complexes. In \cite{harrach_twisted} it was shown (in a much more
  general setting) that a Poisson transform for differential forms has harmonic
  values if and only if it descends to the BGG complex in an appropriate sense, see
  also Section \ref{sec:harmonic} below.

  Building on the work in the second author's thesis \cite{harrach_thesis}, we have
  constructed in the article \cite{cap_harrach_julg} Poisson transforms that descend
  to the BGG complex in the complex case. This was based on an analysis of all
  possible Poisson transforms and rather extensive computations to find appropriate
  combinations. In Section \ref{sec:general} of the current article, we use
  representation theory to devise a construction of appropriate Poisson transforms,
  which applies to all cases of real rank one (and partly in more
  generality). In particular, we prove for all cases of real rank one
    that the values of the resulting transforms are not only harmonic but also
  coclosed and analyze their compatibility with the exterior derivative on hyperbolic
  space and the operators in the BGG complex on the sphere.

  Starting from Section \ref{3}, we specialize to the complex case
    $G=\SU(n+1,1)$. In this case, we study harmonic Poisson transforms in detail,
  analyzing their compatibility with the operations coming from the K\"ahler
  structure on complex hyperbolic space. In particular, we show that the transforms
  in different degrees can be normalized in such a way that one obtains a chain map
  from the BGG complex on the sphere to the de Rham complex of complex hyperbolic
  space. Next, we use the Poincar\'e ball model to realize the sphere as the boundary
  (at infinity) of complex hyperbolic space and study the boundary asymptotics of the
  values of the transforms in this picture. This builds on a careful analysis of the
  restrictions of the transform to isotypical components for the action of the maximal
  compact subgroup. In particular, we prove that on $S^{2n+1}$, for forms of degree
  $\leq n$, the values under the transform admit continuous extensions to the
  boundary. This is in turn used to analyze the composition of $d$ with the transform
  in degree $n$, showing in particular that its values are $L^2$-harmonic forms. This
  leads to Theorem \ref{prop_L_2_general}, which is the main result of the article:
  We show that this composition isomorphically embeds the quotient of the degree
  $n$-part of the BGG complex by the image of the BGG operator landing in this space
  onto a dense subspace of the space of $L^2$-harmonic forms of degree $n+1$ in the
  interior, which is known to be the direct sum of the discrete series
  representations of $SU(n+1,1)$ with trivial infinitesimal character.

  A major motivation for studying Poisson transforms in our setting is the program
  put forward by the third author to use them to prove special cases of the
  Baum-Connes conjecture with coefficients. We briefly describe this here, referring
  to \cite{gjv_Baum_Connes} and \cite{julg_sp} for more details. 

The Baum-Connes conjecture for a locally compact group $G$ aims to describe the
$K$-theory of the reduced $C^*$-algebra of $G$ in geometrical terms. From the
non-commutative geometry point of view, the group $K_*(C^*_r(G))$ is the $K$--theory
of the space of unitary representations of $G$ which appear in the left regular
representation $L^2(G)$. Baum, Connes and Higson have constructed an assembly map
from the $K$-homology of the classifying space for proper actions of $G$ to the group
$K_*(C^*_r(G))$. The conjecture says that this map is an isomorphism. It is important
to note that the case where $G$ is a semisimple Lie group is well understood, since
representation theory is well understood. See the work of A. Wassermann (using
Harish-Chandra theory and the Knapp-Stein intertwining operators) or more recently
the beautiful approach of A. Afgoustidis (using D. Vogan's theory of minimal
$K$-types).

On the other hand, the case of discrete groups is much more problematic, since no
theory of representations is available.  However, an interesting (and in general
still open) problem is the following: given a connected semisimple Lie group $G$, we
would like to prove the Baum-Connes conjecture for \emph{any} discrete subgroup
$\Gamma$ of $G$. This leads to a technical generalization of the conjecture, namely
the Baum-Connes conjecture {\it with coefficients} for a locally compact group $G$
acting by automorphisms on a $C^*$-algebra $A$, which is the analogue of the original
conjecture for the $K$-theory groups $K_*(C^*_r(G,A))$. Using the techniques of
bivariant $K$-theory of G.\ Kasparov, one can show that if a group satisfies BC with
coefficients (i.e. for any action on a $C^*$-algebra), then so do all its closed
subgroups. This applies to our case: if a connected semisimple group $G$ satisfies
the Baum-Connes conjecture with coefficients, then any discrete subgroup $\Gamma$ of
$G$ satisfy the Baum-Connes conjecture.

We are therefore lead to investigate the conjecture {\it with coefficients} for
connected semisimple Lie groups.  Using the geometry of symmetric spaces, Kasparov
proves that the assembly map is always injective, and describes the condition for
surjectivity in terms of a certain element, denoted $\gamma$, of a ring of bivariant
$K$-theory. More precisely, there is an abelian ring $KK_G(\mc,\mc)$ which acts by
endomorphisms on all the $K$-theory groups $K_*(C^*_r(G,A))$. Kasparov constructs an
idempotent element $\gamma\in KK_G(\mc,\mc)$ which has the following
property: the assembly map to $K_*(C^*_r(G,A))$ is surjective if and only if $\gamma$
acts by the identity on $K_*(C^*_r(G,A))$.

Let us be a little more precise. An element of $KK_G(\mc,\mc)$ is a homotopy
class of a $G$-Fredholm module. By $G$-Fredholm module, we mean the following data:
two Hilbert spaces $H_0$ and $H_1$ equipped with unitary representations $\pi_0$ and
$\pi_1$ of $G$, and a Fredholm map $H_0\rightarrow H_1$ such that the commutators
$F\pi_0(g)-\pi_1(g)F$ are compact operators.

The work of many people (J. Fox, P. Haskell, G. Kasparov, Z.Q. Chen and the third
author) in the 1980's and 90's has consisted in constructing an explicit Fredholm
module representing the element $\gamma$ in the case of the groups $SO_0(n+1,1)$ and
$SU(n+1,1)$, and eventually proving the conjecture in these cases. The third author has
generalized such a construction to $Sp(n,1)$, see \cite{julg_sp}, and hopes to
provide a proof of the BC-conjecture with coefficients for that group. The
interesting point is that the unitary representations appearing in the Fredholm
module are sums of a finite number of irreducible representations. Those irreducible
components are either discrete series or principal series representations. More
precisely:

1) The discrete series appearing in the representative of $\gamma$ are the components
of the $L^2$-cohomology of the symmetric space $G/K$, i.e. the Hilbert space of
square integrable harmonic forms on $G/K$. They are absent for $G=SO_0(2n+1,1)$;
there are 2 components in the case of $SO_0(2n,1)$ and $n+2$ components in the cases
of $SU(n+1,1)$ and $Sp(n+1,1)$.

2) The principal series appearing in $\gamma$ are also very specific ones. In the
case of $SO_0(n+1,1)$, they are exactly the differential forms on the boundary sphere
$G/P$ (made unitary after twisting by some density line bundle). But the general case
is more subtle. Let us explain that: let $P=MAN$ be the Iwasawa decomposition of the
parabolic group $P$.  A principal series representation is induced from $P$ to $G$ of
a representation of $P$ which is trivial on the nilpotent group $N$ (equivalently a
representation of $M\times A$). The representation of $G$ on the spaces of
differential forms splits into principal series if and only if the adjoint action of
$N$ on its Lie algebra is trivial, i.e. iff $N$ is abelian. This is the case for
$SO_0(n+1,1)$ since $N$ is isomorphic to $\mr^n$, but not for $SU(n+1,1)$ where
$N$ is the Heisenberg group $\mc^n\oplus i\mr$. The same occurs in the case
of $Sp(n+1,1)$, where $N = \mh^n \oplus \im(\mh)$.  This forces us to replace the bundle of forms on $G/P$ by a
subquotient on which the action of $N$ is trivial: this means exactly to replace the
de Rham complex by the BGG complex.

The last ingredient in the construction of the $G$-Fredholm module is a Poisson map
between the BGG complex on $G/P$ (more precisely the component below the middle) and
the $L^2$-harmonic forms of $G/K$. In the $SO_0(n+1,1)$ case this is exactly the map
from $\Omega^n(G/P)$ to $\Omega^{n+1}(G/K)$ given by composition of P.Y. Gaillard's
Poisson transform in degree $n$ with the de Rham operator on $G/K$ as mentioned
above. In the current paper, we sort out the case of $\SU(n+1,1)$, and make precise
the link with the computations in \cite{julg_kasparov}. As mentioned above,
Section \ref{sec:general} provides several ingredients needed for
attacking the case of $Sp(n+1,1)$.

  \subsection*{Acknowledgements} First and second author supported by the Austrian
  Science Fund (FWF): P 33559-N. This article is based upon work from COST Action
  CaLISTA CA21109 supported by COST (European Cooperation in Science and
  Technology). https://www.cost.eu. We thank the anonymous referee for several 
  comments that helped improving the paper. 

\noindent
\textbf{Competing interests:} The author(s) declare none.

\section{General theory of Poisson transforms}\label{sec:general}
\subsection{Harmonic Poisson transforms}\label{sec:harmonic}
Throughout this section let $\xg$ be a non-compact real semisimple Lie algebra with
Cartan involution $\theta$ and corresponding Cartan decomposition $\xg = \xk \oplus
\xq$ so that $\xq \neq \{0\}$. Choose a maximal abelian subspace $\xa_0 \subset \xq$
and let $\Delta_r \subset \xa_0^*$ be the set of restricted roots of $\xg$ with
respect to $\xa_0$. For $\lambda \in \Delta_r$ let $\xg_{\lambda}$ be its restricted
root space and define $\xm_0$ as the centralizer of $\xa_0$ in $\xk$. Choose subsets
$\Delta_r^+$ and $\Sigma_0$ of positive and simple roots, respectively, and let
$\xn_0 = \bigoplus_{\lambda \in \Delta_r^+} \xg_{\lambda}$. Then $\xp_0 := \xm_0
\oplus \xa_0 \oplus \xn_0$ is a Borel subalgebra of $\xg$. Let $\xp$ be a standard
parabolic subalgebra of $\xg$ with respect to $\xp_0$ and let $\Sigma \subset
\Sigma_0$ be its corresponding subset of simple roots. Let $\xg = \xg_{-k} \oplus
\dotso \oplus \xg_k$ be the $|k|$-grading of $\xg$ induced by $\Sigma$-height. Define
$\xp_+$ and $\xg_-$ to be the sum of all positive and negative grading components,
respectively, then by construction we have $\xp = \xg_0 \oplus \xp_+$. Put $\xm :=
\xk \cap \xp = \xk \cap \xg_0$ as well as $\xa := \xq \cap \xg_0$. Let $G$ be a
connected semisimple Lie group with Lie algebra $\xg$ and assume that the centre of
$G$ is finite. Let $K$ be the maximal compact subgroup with Lie algebra $\xk$, $P =
N_G(\xp)$ be the parabolic subgroup of $G$ with Lie algebra $\xp$ and let $P_+
\subset P$ be its unipotent radical whose Lie algebra coincides with $\xp_+$. We
identify the quotient group $G_0 := P/P_+$ with a subgroup of $G$ so that the Lie
algebra of $G_0$ coincides with $\xg_0$.

We briefly recall the construction of a $G$-equivariant operator mapping between
vector valued differential forms on $G/P$ and $G/K$ as presented in
\cite{harrach_twisted}, see also \cite{harrach_thesis}. The Iwasawa decomposition
yields a diffeomorphism $G/M \cong G/K \times G/P$ of $G$-manifolds, where $M := K
\cap P$. Denoting the natural projections from $G/M$ to $G/K$ and $G/P$ by $\pi_K$
and $\pi_P$, respectively, we obtain an induced decomposition of the tangent bundle
of $G/M$ into the direct sum of $T_K := \ker(T\pi_P)$ and $T_P :=
\ker(T\pi_K)$. Next, for an irreducible $G$-representation $\mv$ and a closed
subgroup $H$ of $G$ we denote by $V_H := G \times_H \mv$ its associated vector bundle
over $G/H$. Now for another irreducible $G$-representation $\mw$ we obtain a
decomposition
$$
 \Lambda^{k}T^*(G/M, L(V_M, W_M)) = \bigoplus_{k_1+k_2=k} \Lambda^{k_1,k_2}T^*(G/M, L(V_M, W_M)),
$$
where the elements in $\Lambda^{k_1,k_2}$ vanish unless exactly $k_1$ of their
entries come from $T_K$ and the remaining $k_2$ entries are from $T_P$.

Assuming that $G/P$ is oriented and $\dim(G/P) = N$, we fix a $G$-invariant
$(0,N)$-form $\vol_M$ on $G/M$, which is unique up to multiples, and normalize it so
that the fibre $\pi_K^{-1}(eK)$ has unit volume. Then each $G$-equivariant bundle map
\begin{equation}\label{bundle-map}
 \phi \colon \Lambda^{0,k} T^*(G/M, V_M) \to \Lambda^{\ell,0} T^*(G/M, W_M) 
\end{equation}
induces, via integration along the fiber, a $G$-equivariant integral operator
\begin{align*}
\Phi &\colon \Omega^k(G/P, V_P) \to \Omega^{\ell}(G/K, W_K), & \Phi(\alpha) =
\fint_{G/P} \phi(\pi_P^*\alpha) \otimes \vol_M,
\end{align*}
called a \emph{Poisson transform}. Equivalently, any such bundle map induces a
$G$-invariant differential form $\varphi$ of bidegree $(\ell, N-k)$ on $G/M$ with
values in $L(V_M, W_M)$ via the relation
\begin{equation}\label{eq:kernel-homom}
 \varphi \wedge \pi_P^*\alpha = \phi(\pi_P^*\alpha) \otimes \vol_M.
\end{equation}

For each closed subgroup $H$ of $G$ the homogeneous bundle $V_H = G \times_H \mv$
carries a canonical $G$-invariant flat connection $\nabla^{V_H}$, which induces a
covariant exterior derivative $d^{V_H}$ on $\Omega^*(G/H, V_H)$ in the usual
way. Additionally, we can find an inner product on $\mv$ so that
$\langle X v_1, v_2 \rangle = \langle v_1, \theta(X) v_2\rangle$ for all $X \in \xg$
and $v_1$, $v_2 \in \mv$. By construction this inner product is $K$-invariant and
thus yields a $G$-invariant bundle metric on $V_K \to G/K$. Next, we define
$\nabla^{\theta}$ as the dual connection to $\nabla^{V_K}$ with respect to this
bundle metric, which is again flat. Moreover, we obtain an induced covariant exterior
derivative $d^{\theta}$ on $\Omega^*(G/K,V_K)$ and we define $\delta^{V_K}$ as the
formal adjoint of $d^{\theta}$ with respect to the induced $L^2$-norm. Finally,
define the Laplace operator $\Delta^{V_K}$ via
$\Delta^{V_K} := \delta^{V_K}d^{V_K} + d^{V_K}\delta^{V_K}$. In
\cite{harrach_twisted}*{Theorem 6} it was proven

\begin{prop}\label{prop_Poisson_differential_operators}
For each Poisson transform $\Phi \colon \Omega^k(G/P, V_P)\to \Omega^\ell(G/K, W_K)$ the compositions $\Phi \circ d^{V_P}$, $d^{W_K} \circ \Phi$, $\delta^{W_K} \circ \Phi$ and $\Delta^{W_K} \circ \Phi$ are again Poisson transforms.
\end{prop}

In the next step we recall the BGG-calculus on $G/P$ and afterwards consider Poisson
transforms which are naturally linked to BGG-complexes. Consider the standard complex
$(C^*(\xg_-, \mv), \frak d)$ for computing Lie algebra cohomology of the nilpotent
Lie algebra $\xg_-$ with values in a $\xg$-representation $\mv$. Each chain space is
naturally a $G_0$-module and the differentials are $G_0$-equivariant. On the other
hand, under the identification $\xg_-^* \cong \xp_+$ the Kostant codifferential 
$\frak d^* \colon C^k(\xg_-, \mv) \to C^{k-1}(\xg_-, \mv)$ introduced originally in
\cite{kostant} is defined via 
\begin{align}\label{eq:Kostant}
\begin{split}
 \frak d^*(Z_1 \wedge \dotso \wedge Z_k \otimes v) &:=\ \sum_{i=1}^k (-1)^i Z_1\wedge \dotso \wedge \widehat{Z_i} \wedge\dotso \wedge Z_k \otimes Z_i \cdot v \\
 &+ \sum_{i<j} (-1)^{i+j} [Z_i, Z_j] \wedge Z_0 \wedge \dotso \wedge \widehat{Z_i} \wedge \dotso \wedge \widehat{Z_j} \wedge \dotso \wedge Z_k \otimes v
\end{split}
\end{align}
for all $Z_i \in \xp_+$ and $v \in \mv$, where the hat denotes omission. It is
immediate that this map is $P$-equivariant and squares to zero. Therefore, the
induced homology spaces $H_k := H_k(\xp_+, \mv) := \ker(\frak d^*)/\im(\frak d^*)$ 
are naturally $P$-modules which are completely reducible,
c.f. \cite{cap_slovak}*{Corollary 3.3.1}. Let $\Box = \frak d\frak d^* +
\frak d^*\frak d$ be the Kostant Laplace on $C^k(\xg_-, \mv)$, which is
$G_0$-equivariant. We obtain a $G_0$-invariant decomposition
\begin{equation}\label{eq:decomp_Laplace}
 C^k(\xg_-, \mv) = \im(\frak d) \oplus \ker(\Box) \oplus \im(\frak d^*),
\end{equation}
and the sum of the first two and the last two summands coincide with $\ker(\frak d)$
and $\ker(\frak d^*)$, respectively. In particular, the homology spaces $H_k$ are
isomorphic to $\ker(\Box)$ as $G_0$-modules. 

Denote the $G$-equivariant bundle map on $\Lambda^* T^*(G/P)\otimes V_P$ induced by
$\frak d^*$ by the same symbol. Moreover, let $\ch_k = \ch_k(G/P, V_P)$ be the
homogeneous vector bundle over $G/P$ associated to $H_k$. Since the operator $\Box$
is $G_0$-equivariant we obtain an induced tensorial operator on the associated graded
of $\Lambda^k T^*(G/P)\otimes V_P$ and the decomposition \eqref{eq:decomp_Laplace}
globalizes only to this space. To remedy this fact let $\Box^R := d^{V_P}\frak d^* +
\frak d^* d^{V_P}$ be the curved Box operator, then from
\cite{calderbank_diemer}*{Theorem 5.2} we know that $\Box^R$ is invertible on
sections of the subbundle induced by $\im(\frak d^*)$ and the inverse $Q^R$ can be
written as a universal polynomial in $\frak d^*d^{V_P}$ by
\cite{cap_relative}*{Theorem 3.9}. Using this it can be shown that
\begin{equation}\label{eq:decomp_Laplace_curved}
 \Omega^k(G/P, V_P) = \im(d^{V_P}\frak d^*) \oplus \ker(\Box^R) \oplus
 \Gamma(\im(\frak d^*)),
\end{equation}
(c.f. \cite{harrach_twisted}*{p. 24}) and the projection $\Pi^R \colon \Omega^k(G/P,
V_P) \to \ker(\Box^R)$ onto the middle summand is given by
\begin{equation}\label{eq:pi_R_formula}
  \Pi^R(\alpha) := \alpha - Q^R\frak d^*d^{V_P}\alpha - d^{V_P}Q^R\frak d^*\alpha
\end{equation}
for all $\alpha \in \Omega^k(G/P, V_P)$. From the definition of $Q^R$,
  one easily concludes that \eqref{eq:pi_R_formula} implies $\frak d^* \circ \Pi^R =
  0$. In particular, $\Pi^R$ factors to an operator $L \colon \Gamma(\ch_k) \to
\Gamma(\ker(\frak d^*))$ which is a $G$-equivariant splitting operator for the
canonical projection $\pi_{\ch}:\Gamma(\ker(\frak d^*))\to\Gamma(\ch_k)$ and it is
characterized by $\frak d^*d^{V_P}L = 0$. We define the \emph{$k$-th
  BGG-operator} $D_k \colon \Gamma(\ch_{k-1}) \to \Gamma(\ch_k)$ as the composition
$D_k := \pi_{\ch} d^{V_P} L$. Since the connection $\nabla^{V_P}$ on
$V_P$ is flat it follows that the BGG-sequence
$$
 \xymatrix{ 0 \ar[r] & \Gamma(\ch_0) \ar[r]^-{D_1} & \Gamma(\ch_1) \ar[r]^-{D_2} & \dotso \ar[r]^-{D_{N}} & \Gamma(\ch_N) \ar[r] & 0}
 $$
is a complex which computes the twisted deRham cohomology of $G/P$ with values in
 the vector bundle $V_P$, c.f. \cite{cap_slovak_soucek}*{Corollary 4.15}.

Our aim is to construct Poisson transforms which satisfy $\Phi\circ \frak d^* =
0$ as well as $\Phi \circ d^{V_P} \circ \frak d^* = 0$. By the first property $\Phi$
naturally factors to an operator $\underline{\Phi} \colon \Gamma(\ch_k) \to
\Omega^\ell(G/K, W_K)$ defined by $\underline{\Phi}(\sigma) := \Phi(\alpha)$, where
$\alpha \in \Gamma(\ker(\frak d^*))$ is any representative of $\sigma \in
\Gamma(\ch_k)$. The second property yields compatibility with the $k$-th
BGG-operator, i.e. for all $\tau \in \Gamma(\ch_{k-1})$ the value
$\underline{\Phi}(D_k\tau)$ coincides with $\Phi(d^{V_P}\beta)$ for any
representative $\beta \in \Gamma(\ker(\frak d^*))$ of $\tau$,
c.f. \cite{harrach_twisted}*{Section 3.1}.

In order to relate the BGG-complex to the geometry of $G/K$ recall that for each
homogeneous vector bundle $E \to G/H$ we have a natural $G$-action on
$\Gamma(E)$. Infinitesimally, one obtains an action of the universal enveloping
algebra $\cU(\xg)$ on sections of $E$, so in particular the Casimir element $\cc \in
\cU(\xg)$ induces a $G$-invariant differential operator $\Gamma(E) \to \Gamma(E)$
which we denote by the same symbol.

\begin{thm}\label{thm_Poisson_BGG}
 Let $\Phi \colon \Omega^k(G/P, V_P) \to \Omega^{\ell}(G/K, W_K)$ be a Poisson transform.
 \begin{enumerate}[(i)]
  \item Let $c_{\mv}$ and $c_{\mw}$ be the Casimir eigenvalues of $\mv$ and $\mw$, respectively. Then 
  \begin{align*}
   2\Phi(\Box^R\alpha) + \Delta^{W_K}\Phi(\alpha) = (c_{\mw} - c_{\mv}) \Phi(\alpha).
  \end{align*}
In particular, if $\Phi \circ \Box^R = 0$, then for any $\alpha \in
\Omega^k(G/P,V_P)$ the element $\Phi(\alpha)$ is an eigenform of $\Delta^{W_K}$ with
eigenvalue $(c_{\mw} - c_{\mv})$, or equivalently, an eigenform of $\cc$ with
eigenvalue $c_{\mv}$.
\item The following are equivalent: \\
\begin{enumerate*}[label = (\alph*), itemjoin = {\quad}]
\item $\Phi \circ \Box^R = 0$,
\item $\Phi \circ \frak d^* = 0$ and $\Phi \circ d^{V_P} \circ \frak d^* = 0$,
\item $\Phi = \Phi \circ \Pi^R$.
\end{enumerate*}

If these are satisfied, $\Phi$ is uniquely determined by its induced operator
$\underline{\Phi}$ on $\Gamma(\ch_k)$.
\end{enumerate}
 \end{thm}
\begin{proof}
 (i) From \cite{cap_soucek_casimir}*{Corollary 1} we know that on $\Omega^k(G/P, V_P)$ we have $\cc = 2 \Box^R + c_{\mv}$. On the other hand, in \cite{matsushima_murakami}*{p. 385, formula (6.9)} it was proven that $\cc = -\Delta^W + c_{\mw}$ on $\Omega^*(G/K, W_K)$. Since Poisson transforms are $G$-equivariant they commute with $\cc$, implying the claim.
 
(ii) The decomposition \eqref{eq:decomp_Laplace_curved} immediately implies equivalence of $(b)$ and $(c)$, whereas equivalence of $(a)$ and $(b)$ was shown in \cite{harrach_twisted}*{Theorem 4}. Finally, for the uniqueness result note that the splitting operator is an isomorphism $\Gamma(\ch_k) \to \Gamma(\ker(\Box^R))$ and satisfies $\underline{\Phi} = \Phi \circ L$. 
 \end{proof}

Part (i) of the previous Theorem shows that for $\mw = \mv$ the condition $\Phi \circ \Box^R = 0$ is equivalent to $\Delta^{V_K} \circ \Phi = 0$. In light of this we have the following definition.
\begin{defin}
 We say that a Poisson transform $\Phi \colon \Omega^k(G/P,V_P) \to \Omega^{\ell}(G/K, W_K)$ is \emph{harmonic} if it satisfies one of the three equivalent properties of part (ii) in Theorem \ref{thm_Poisson_BGG}.
\end{defin}

\subsection{BGG-compatibility and equivariant maps.}\label{sec_BGG_equiv}
A crucial feature of the approach to Poisson transforms introduced above
  is that in view of homogeneity of $G/M$, $G$-equivariant bundle maps as
  in \eqref{bundle-map} are in bijective correspondence with $M$-equivariant maps
  between the representations inducing the bundles. We next describe harmonic Poisson
  transforms in this language.

  We start with some general observations. Since $M=K\cap P$, any representation
  $\mv$ of $K$ restricts to a representation of $M$ and hence gives rise to a
  homogeneous bundle $G\x_M\mv\to G/M$, which can be identified with the pullback
  along $\pi_K$ of the homogeneous bundle $G\x_K\mv\to G/K$. Natural bundle maps
  between homogeneous bundles over $G/K$ come from $K$-equivariant maps between the
  inducing representations and hence lift to natural bundle maps between the pullback
  bundles. This also applies to the projection to $G/P$ but in this case, more is
  true. Since $M\subset G_0\subset P$, also representations of $G_0$ give rise to
  natural bundles over $G/M$ and $G_0$-equivariant maps induce natural bundle maps
  between such homogeneous bundles. In the case of representations of $G_0$ which are
  not restrictions of representations of $P$, one has to be careful with the pullback
  interpretation, though. This is available for associated bundles to the principal
  $G_0$-bundle $G/P_+\to G/P$ which can be interpreted as homogeneous bundles
  corresponding to completely reducible representations of $P$.

  The $M$-representation inducing $T(G/M)$ is $\xg/\xm$. The subbundle $T_K$
  corresponds to $(\xg/\xm)^{1,0} := \xp/\xm$ which is isomorphic to the restriction
  to $M$ of the $K$-representation $\frak g/\frak k$. Likewise, $T_P$ corresponds to
  $(\xg/\xm)^{0,1} := \xk/\xm$ which is isomorphic to the restriction of the
  $P$-representation $\xg/\xp$. For the induced decomposition of the dual, there is a
  more natural point of view as $(\xg/\xm)^*=(\xg/\xk)^*\oplus(\xg/\xp)^*$. Viewing
  $(\xg/\xm)^*$ as the annihilator of $\xm$ in $\xg^*$, the summands are naturally
  contained as annihilators of bigger subalgebras. Hence
  $\Lambda^{k_1,k_2} (\xg/\xm)^* := \Lambda^{k_1} (\xg/\xk)^* \otimes \Lambda^{k_2}
  (\xg/\xp)^*$ is the $M$-module underlying $\Lambda^{k_1,k_2} T^*(G/M)$. The above
  observations on restrictions of $G_0$-representations to $M$ apply to the chain
  spaces $C^k(\xg_-, \mv)$ and the equivariant maps $\frak d$, $\frak d^*$ and
  $\Box$, for which we denote the induced natural bundle maps by $\frak d_P$,
  $\frak d_P^*$ and $\Box_P$. The splitting \eqref{eq:decomp_Laplace} then reads as
\begin{align}\label{eq:decomp_Laplace_G_M}
 \Lambda^{0,k}T^*(G/M)\otimes V_M = \im(\frak d_P) \oplus \ker(\Box_P) \oplus
 \im(\frak d_P^*),
\end{align}
and the sum of the first and last two modules coincide with $\ker(\frak d_P)$ and
$\ker(\frak d_P^*)$, respectively. Finally, $H_k = H_k(\xp_+,
\mv)\cong_M\ker(\Box_P)$, so over $G/M$ we can identify $\Gamma(\ch_k)$ with
$\Gamma(\ker(\Box_P))$ and its elements are described by $M$-equivariant maps $G\to
\ker(\Box_P)$.

\begin{thm}\label{thm_uniqueness}
 Let $\Phi \colon \Omega^k(G/P, V_P) \to \Omega^{\ell}(G/K, W_K)$ be a Poisson transform.
 \begin{enumerate}[(i)]
 \item The compositions $\Phi \circ \frak d^*$ and $\Phi \circ \Pi^R$ are again Poisson transforms.
 
 \item If $\Phi \circ \frak d^* = 0$ then $\Psi := \Phi \circ \Pi^R$ is a harmonic
   Poisson transform which induces the same transform on $\Gamma(\ch_k)$ as
   $\Phi$. Moreover, if $\sigma \in \Gamma(\ch_k)$ corresponds to $\underline{\sigma}
   \in \Gamma(\ker(\Box_P))$, then
 \begin{align*}
  \underline{\Phi}(\sigma) = \fint_{G/P} \psi(\underline{\sigma}) \otimes \vol_M,
 \end{align*}
where $\psi$ is the $M$-equivariant map underlying $\Psi$. Therefore, harmonic Poisson transforms are in bijective correspondence with $\Hom_M(\ker(\Box_P), \Lambda^{\ell,0}(\xg/\xm)^* \otimes \mw)$.
  \end{enumerate}
\end{thm}
\begin{proof}
 (i) By construction $\frak d_P^*\pi_P^*\alpha$ coincides with
    $\pi_P^*\frak d^*\alpha$ for all $\alpha \in \Omega^k(G/P, V_P)$, and thus $\Phi
    \circ \frak d^*$ is a Poisson transform. Since $\Phi \circ d^{V_P}$ is a Poisson
    transform due to Proposition \ref{prop_Poisson_differential_operators}, the same
    is true for $\Phi \circ \Pi^R$ as $\Pi^R$ is a polynomial in $\frak d^*d^{V_P}$.
  
  (ii) From Theorem \ref{thm_Poisson_BGG}(ii), we know that the operator $\Psi = \Phi
    \circ \Pi^R$ is harmonic. Moreover, given $\sigma \in \Gamma(\ch_k)$ with
    representative $\alpha \in \Gamma(\ker(\frak d^*))$, then $\Pi^R(\alpha)$ is
    also a representative of $\sigma$. This implies that the maps induced by $\Phi$
    and $\Psi$ on $\Gamma(\ch_k)$ coincide, and uniqueness follows directly from
    Theorem \ref{thm_Poisson_BGG}(ii). Next, $\Psi \circ \frak d^* = 0$ implies that
    $\psi$ is trivial on $\im(\frak d_P^*)$, so from the decomposition
    \eqref{eq:decomp_Laplace_G_M} it is immediate that $\psi(\pi_P^*\alpha)$
    coincides with $\psi(\underline{\sigma})$. This also shows that $\Psi$ is fully
    determined by the restriction of $\psi$ to $\ker(\Box_P)$.  
\end{proof}

Recall that for any $G$-representation $\mv$, the homology space
$H_k = H_k(\xp_+, \mv)$ decomposes into a direct sum of irreducible
$P$-representations $\mE_i^k$ and we denote by $\pi^k_i$ the projection onto the
component $\mE_i^k$. In the same way, the $M$-representation $\ker(\Box_P)$
decomposes into irreducible components that we also denote by
$\mE_i^k$. Correspondingly $\ch_k = G \times_P H_k=\oplus_iE_i^k$, where
$E_i^k := G \times_P \mE_i^k$, and the projection onto one summand is a
$G$-equivariant bundle map. We also denote this bundle map and the induced tensorial
operator on sections by $\pi_i^k$. This leads to a splitting of the $k$th
BGG-operator into components $D^{j,i}_k$ defined by
$D^{j,i}_k := \pi_i^k \circ D_k\circ \pi_j^{k-1}$.

  Now suppose that we have given a harmonic Poisson transform
  $\Phi:\Omega^k(G/P,V_P)\to\Omega^\ell(G/K,W_K)$ induced by
  $\phi\in\Hom_M(\ker(\Box_P), \Lambda^{\ell,0} (\xg/\xm)^*\otimes\mw)$. Then we
  define $\phi_i:=\phi\o\pi_i^k$ and this induces a harmonic Poisson transform
  $\Phi_i$ acting between the same spaces as $\Phi$. Part (ii) of Theorem
  \ref{thm_uniqueness} readily shows that for the induced transform we get
  $\underline{\Phi}_i=\underline{\Phi}\o\pi_i^k$. Using this, we can now show that
  harmonic Poisson transforms are well behaved with respect to compositions with
  components of BGG operators.

  \begin{prop}\label{prop_poisson_on_irreducible}
 Let $\Phi \colon \Omega^k(G/P,V_P) \to \Omega^{\ell}(G/K, W_K)$ be a harmonic
 Poisson transform and let $\Phi_i$ be the component corresponding to $E_i^k \subset
 \ch_k$. Consider the $k$th BGG operator $D_k$ and  a component $D^{j,i}_k$ with
 values in $\Gamma(E_i^k)$. 
 \begin{enumerate}[(i)]
 \item The transform $\underline{\Phi}\o D$ on $\Gamma(\ch_{k-1})$ is induced by a
   harmonic Poisson transform defined on $\Omega^{k-1}(G/P,V_P)$. 
  \item Let $\Psi$ be a harmonic Poisson transform which in the setting of (i)
    induces $\underline{\Phi}_i\o D$. Then the harmonic Poisson transform $\Psi_j$
    induces $\underline{\Phi}\o D^{i,j}_k:\Gamma(\ch_{k-1})\to\Omega^{\ell}(G/K, W_K)$.
 \end{enumerate}
  \end{prop}
  \begin{proof}
(i) By definition of the BGG operators, for any $\al\in\Gamma(\ch_{k-1})$, the
       differential form $d^{V_P}(L(\al))$ is a representative for
       $D(\al)\in\Gamma(\ch_k)$, and hence $\Phi\o d^{V_P}\o \Pi^R$ induces
       $\underline{\Phi}\o D$ on homology bundles. Now $\Phi\o d^{V_P}$ is a Poisson
       transform by Proposition \ref{prop_Poisson_differential_operators} and since
    $\Phi$ is harmonic, we get $(\Phi\o d^{V_P})\o\frak d^*=0$. By part (i) of
       Theorem \ref{thm_uniqueness} this implies that $\Phi\o d^{V_P}\o \Pi^R$ is a
       harmonic Poisson transform.

(ii) By construction, we have 
$\underline{\Psi}(\be)=\underline{\Phi}(\pi^k_iD(\be))$ and
$\underline{\Psi}_j(\be)=\underline{\Psi}(\pi_j^{k-1}(\be))$ which implies the claim.
  \end{proof}

  \subsection{A specific family of transforms}\label{sec:our-Poisson}
  We now specialize to a family of transforms that will be used in what follows. We
  do this only in the case that $\mv$ and $\mw$ are $\mr$ or $\mc$ and that $\xg$ is
  of real rank one. This implies that $\xp\subset\xg$ is a minimal parabolic
  subalgebra and that $\xa=\xq\cap\xg_0$ is a one dimensional abelian subalgebra
  spanned by the grading element $E$. This element is characterized by $[E, X] = jX$
  for all $X \in \xg_j$, which implies that $M$ acts trivially on $\xa$.

  It is well known that $\xp_+$ is the annihilator of $\xp$ under the Killing form
  $B$, so one obtains an isomorphism $(\xg/\xp)^*\cong\xp_+$. This is $P$-equivariant
  and hence $M$-equivariant and thus can also be viewed as an isomorphism
  $\Lambda^{0,1}(\xg/\xm)^*\cong_M\xp_+$. Similarly, $\xq$ is the annihilator of
  $\xk$ under the Killing form, so we get an isomorphism $(\xg/\xk)^*\cong\xq$. This
  is $K$-equivariant and hence $M$-equivariant, so
  $\Lambda^{1,0}(\xg/\xm)^*\cong_M\xq$. But now $\xq\cap(\xg_-\oplus\xp_+)$ is an
  $M$-invariant complement to $\xa$ in $\xq$ and any element in there can be uniquely
  written as $\theta(Z)-Z$ for some $Z\in\xp_+$. Thus we conclude that
  $\Lambda^{1,0}(\xg/\xm)^*\cong \xa\oplus\xp_+$ as a representation of
  $M$. Otherwise put, sending $Z\in\xp_+$ to the map $Y\mapsto B(Y,\theta(Z)-Z)$
  defines an $M$-equivariant injection $\xp_+\to \Lambda^{1,0}(\xg/\xm)^*$. Together
  with the above, this defines an $M$-equivariant injection
  $\iota_1:\Lambda^{0,1}(\xg/\xm)^*\to \Lambda^{1,0}(\xg/\xm)^*$. Taking exterior
  powers, we obtain $M$-equivariant maps $\iota_k:\Lambda^{0,k}(\xg/\xm)^*\to
  \Lambda^{k,0}(\xg/\xm)^*$.

  It is easy to make this map explicit. An element
  $\alpha\in\Lambda^{0,k}(\xg/\xm)^*$ is a $k$-linear alternating map $\xg^k\to\mk$,
  where $\mk=\mr$ or $\mc$, which vanishes if one of its entries is from
  $\xp$. Invariance of the Killing form readily implies that for $Z\in\xp_+$ and
  $Y\in\xg$ we get $B(Y,\theta(Z)-Z)=B(\theta(Y)-Y,Z)$. Passing to exterior powers,
  we readily conclude that
  \begin{equation}\label{eq:ik-formula}
    \iota_k(\alpha)(Y_1,\dots,Y_k)=\alpha(\theta(Y_1)-Y_1,\dots,\theta(Y_k)-Y_k). 
  \end{equation}
  Observe that, as expected, this vanishes if one of its entries is in $\xk$. As an
  element of $\Lambda^{k,0}(\xg/\xm)^*$, $\iota_k(\alpha)$ is uniquely determined by
  its value on entries from $\xp$ and for $Y_i\in\xp$, the formula simplifies to 
  \begin{equation}\label{eq:ik-formula-p}
  \iota_k(\alpha)(Y_1,\dots,Y_k)=\alpha(\theta(Y_1),\dots,\theta(Y_k)),
  \end{equation}
  which explains our choice of sign for the map $\iota_1$.
  
  \begin{defin}\label{def-our-Poisson}
    Take $\mk=\mr$ or $\mc$ and let $\mE^k_i\subset\ker(\square_P)\cong H_k(\xp_+,\mk)$
    be an irreducible component. Then we define a map $\phi^k_i :=
    \iota_k\o\pi^k_i\in\Hom_M(\ker(\square_P),\Lambda^{k,0}(\xg/\xm)^*)$. By
    $\Phi^k_i:\Omega^k(G/P,\mk)\to\Omega^k(G/K,\mk)$ we denote the corresponding
    Poisson transform, which is harmonic by Theorem \ref{thm_uniqueness}.
  \end{defin}

  From Section \ref{sec_BGG_equiv} we know that the transform $\Phi^k_i$ defines a
  map $\Gamma (E^k_i)\to\Omega^k(G/K,\mk)$ and we next compute $d\o \Phi^k_i$ and
  $\delta\o \Phi^k_i$. By Theorem 6 of \cite{harrach_twisted} these compositions are
  Poisson transforms induced by the kernels $d_K\varphi^k_i$ and
  $\delta_K\varphi^k_i$, respectively, where $\varphi^k_i$ is the kernel of
  $\Phi^k_i$. Here $d_K$ is the component of the exterior derivative which raises 
  the bidegree by $(1,0)$ and $\delta_K$ corresponds to the codifferential on $G/K$,
  see \cite{harrach_twisted} for details. Since $\varphi^k_i$ is a $G$-invariant
  differential form, these operations can be computed purely algebraically from the
  inducing $M$-invariant elements.

  To do this, we use the Lie algebra cohomology differential $\frak
  d_{\xg}:\Lambda^k\xg^*\to\Lambda^{k+1}\xg^*$ for $\xg$ with values in the trivial
  representation $\mk$. For $\alpha\in\Lambda^k\xg^*$ this is explicitly given by
  \begin{equation}\label{eq:delg}
\frak
d_{\xg}\al(Y_0,\dots,Y_k)=\textstyle\sum_{i<j}(-1)^{i+j}\alpha([Y_i,Y_j],Y_0,\dots,\widehat{Y_i}, \dots, \widehat{Y_j},\dots, Y_k),
  \end{equation}
  where the hats denote omission. Now we first prove a technical lemma.

  \begin{lemma}\label{lem:delg}
    Suppose that $\alpha\in\Lambda^{0,k}(\xg/\xm)^*$ satisfies $\frak d_P\alpha=0$
    and $\beta\in\Lambda^{0,\ell}(\xg/\xm)^*$ satisfies $E\cdot\beta=\nu\be$ for
    $\nu\in\mz$, and consider $\gamma:=\frak
    d_{\xg}(i_k(\alpha)\wedge\beta)\in\Lambda^{k+\ell+1}\xg^*$. Take elements
    $Y_0,\dots Y_{k+\ell}\in\xp_+$ such that $Y_j\in\xg_{s_j}$ for $j>k$ and put
    $\tilde Y_i:=\theta(Y_i)+Y_i$.

    Then $\gamma(Y_0,\dots,Y_k,\tilde Y_{k+1},\dots,\tilde Y_{k+\ell})=0$ provided
    that $s_{j+1}+\dots+s_{k+\ell}\leq\nu$. Moreover, we get
    \begin{equation}\label{eq:tech-delg}
        \gamma(E,Y_1,\dots,Y_k,\tilde Y_{k+1},\dots ,\tilde
        Y_{k+\ell})=\big(-i_k(E\cdot\alpha)\wedge\beta+\nu i_k(\alpha)\wedge
        \beta\big) (Y_1,\dots ,\tilde Y_{k+\ell})
    \end{equation}
  \end{lemma}
  \begin{proof}
   Computing $\gamma(Y_0,\dots,Y_k,\tilde Y_{k+1},\dots,\tilde Y_{k+\ell})$ via
   formula \eqref{eq:delg}, we have to distinguish cases according to which type of
   elements go into the Lie bracket. We know that $i_k(\alpha)$ vanishes if one of
   its entries is from $\xk$ while $\beta$ vanishes if one of its entries is from
   $\xp$. Now $[\tilde Y_i,\tilde Y_j]\in\xk$ so in the corresponding terms we feed
   $k+1$ entries from $\xp_+$ and $\ell-1$ entries from $\xk$ into
   $i_k(\alpha)\wedge\beta$, which automatically gives zero. Similarly,
   $[Y_i,Y_j]\in\xp$ and hence for the corresponding summands, all the $\tilde Y$'s
   have to be inserted into $\beta$. Thus for this summand, 
   we have to consider
   $$
   \left(\textstyle\sum_{i<j\leq k}(-1)^{i+j}i_k(\alpha)([Y_i,Y_j],Y_0,\dots, \hat i,\dots,
   \hat j,\dots, Y_k)\right)\beta(\tilde Y_{k+1},\dots,\tilde Y_{k+\ell}). 
   $$

   Using formula \eqref{eq:ik-formula-p} and that
   $\theta([Y_i,Y_j])=[\theta(Y_i),\theta(Y_j)]$ we conclude that these terms add up
   to zero since $\frak d_P\alpha=0$. Thus we are left with the summands that contain
   $[Y_i,\tilde Y_j]$, for which there are $k$ more entries from $\xp_+$. Hence all
   the remaining $Y$'s have to go into $i_k(\alpha)$ while the bracket term and the
   remaining $\tilde Y$'s have to be inserted into $\beta$. Now $[Y_i,\tilde
     Y_j]=[Y_i,\theta(Y_j)]+[Y_i,Y_j]$ and the second summand lies in $\xp_+$ and 
   hence inserts trivially into $\beta$. For the same reason we can also replace the
   remaining entries of $\beta$ by $\theta(Y_r)$ instead of $\tilde Y_r$.

   By assumption $\theta(Y_r)\in\xg_{-s_r}$ and since $Y_i\in\xp_+$ the bracket
   $[Y_i,\theta(Y_j)]$ is a sum of elements of $\xg_i$ with $i>-s_j$. Hence the
     degrees of the elements inserted into $\beta$ add up to a number which is
     $>-s_{k+1}-\dots -s_{k+\ell}\geq -\nu$. But since we have assume that
     $E\cdot\beta=\nu\beta$, $\beta$ can only be non-zero on a family of elements if
     these degree add up exactly to $-\nu$. Hence these terms do not contribute
     either, and the proof of the first part is complete.

   For the second part, we first observe that $E$ inserts trivially into
   $i_k(\alpha)\wedge\beta$, so in the expansion according to \eqref{eq:delg}, we
   only have to consider the terms in which $E$ gets inserted into the Lie bracket.
   Moreover, $[E,Y_i]\in\xp$ so this bracket has to go into $i_k(\alpha)$ together
   with the other $Y$'s, while all the $\tilde Y$'s have to go into $\beta$. So we
   have to consider
   $$
   -\textstyle\sum_{i=1}^ki_k(\alpha)(Y_1,\dots,[E,Y_i],\dots,Y_k)=
   -\textstyle\sum_{i=1}^k\alpha(\theta(Y_1),\dots,\theta([E,Y_i]),\dots,\theta(Y_k)). 
   $$
   Using that $\theta([E,Y_i])=-[E,\theta(Y_i)]$ we see that this sum equals
   $-(E\cdot\alpha)(\theta(Y_1),\dots,\theta(Y_k))$. So the first $k$ summands in
   \eqref{eq:delg} add up to $-(i_k(E\cdot\alpha)\wedge\beta)(Y_1,\dots,\tilde
   Y_{k+\ell})$. For the remaining terms, $[E,\tilde Y_j]=[E,\theta(Y_j)]+[E,Y_j]$
   has to be inserted into $\beta$, so the second summand does not contribute, while
   all the $Y_i$ have to go into $i_k(\alpha)$. Thus we have to consider the sum
   $$
   -\textstyle\sum_{j=k+1}^{k+\ell}\beta(\theta(Y_{k+1}),\dots,[E,\theta(Y_j)],\dots,
   \theta(Y_{k+\ell}))=(E\cdot\beta)(\tilde Y_{k+1},\dots,\tilde Y_{k+\ell}), 
   $$ so these summands provide the rest of the claimed formula.
  \end{proof}
  
  Viewed as an element of $\xa\subset\Lambda^{1,0}(\xg/\xm)$, the grading element $E$
  defines an $M$-invariant element $E^*\in\Lambda^{1,0}(\xg/\xm)^*$ and hence an
  element of $\Omega^{1,0}(G/M)^G$ that we denote by the same symbol. Using this, we
  can now formulate

  \begin{thm}\label{thm_PT_rank_1}
    Let $\Phi_i^k$ be one of the transforms from Definition \ref{def-our-Poisson}
    induced by the homomorphism $\phi_i^k$. Then we have
 \begin{enumerate}[(i)]
  \item $\delta \circ \underline{\Phi}_i^k=0$
  \item The homomorphism underlying $d \circ \underline{\Phi}_i^k$ is given by
    $2(\rho - \mu_i) E^* \wedge \phi_i^k$, where $2\rho$ and $\mu_i \in \mz$ are the
    eigenvalues of $E$ on $\Lambda^N \xp_+$ and on $\mE^k_i$, respectively.
 \end{enumerate}
\end{thm}
  \begin{proof}
   As we have observed above already, it suffices to show that the kernel
   $\varphi_i^k$ of $\Phi_i^k$ satisfies $d_K\varphi_i^k= 2(\rho - \mu_i) E^* \wedge
   \varphi_i^k$ and $\delta_K\varphi_i^k=0$ to complete the proof. Now it is well
   known (see e.g.\ \cite{onishchick}) that the map $\frak d_{\xg}$ defined in
   formula \eqref{eq:delg} descends to a well defined map
   $(\Lambda^r(\xg/\xm)^*)^M\to (\Lambda^{r+1}(\xg/\xm)^*)^M$, which encodes the
   exterior derivative $\Omega^r(G/M)^G\to\Omega^{r+1}(G/M)^G$.

   We next describe the $M$-invariant element $\psi_i^k\in\Lambda^{k,N-k}(\xg/\xm)^*$
   corresponding to $\phi_i^k$ explicitly. Recall that $(X,Y)\mapsto B(X,\theta(Y))$
   defines a positive definite, $K$-invariant inner product on $\xg$ which restricts
   to $M$-invariant inner products on $\xg_-$ and $\xp_+$ which in turn give rise to
   inner products on the chain spaces $\Lambda^{0,*}(\xg/\xm)^*$. The
   decomposition of the chain spaces corresponding to \eqref{eq:decomp_Laplace_G_M}
   are orthogonal with respect to this inner product, which we denote by
   $\langle\ ,\ \rangle$. (Indeed, Kostant originally defined $\frak d^*$ as an
   adjoint of $\frak d$ with respect to this inner product.) Now we take an
   orthonormal basis $\{\beta_s\}$ for $\Lambda^{0,k}(\xg/\xm)^*$ such that the first $L$
   elements form a basis for $\mE_i^k\subset\ker(\square_P)$, the next elements
   extend this to a basis of $\ker(\square_P)$ and the remaining elements are either
   from $\im(\frak d_P)$ or from $\im(\frak d^*_P)$. Then there are unique
   elements $\beta_s^\#\in \Lambda^{0,N-k}(\xg/\xm)^*$ such that
   $\beta_s\wedge\beta^\#_t=\delta_{s,t}\vol_M$, where we denote by $\vol_M$ also the
   element of $\Lambda^{0,N}(\xg/\xm)^*$ corresponding to the $G$-invariant form with
   the same name.

   For $\alpha\in \Lambda^{0,k}(\xg/\xm)^*$ we conclude that
   $\sum_{i=1}^L\langle\alpha,\beta_s\rangle\beta_s$ is the projection of $\alpha$
   onto $\mE_i^k$ used in Definition \ref{def-our-Poisson}. This immediately implies
   that $\phi_i^k(\alpha)=\sum_{s=1}^L\langle\alpha,\beta_s\rangle i_k(\beta_s)$ and
   in turn from \eqref{eq:kernel-homom} we conclude that
   $\psi_i^k=\sum_{s=1}^Li_k(\beta_s)\wedge\beta^\#_s$. Thus $d_K\varphi_i^k$
   corresponds to the component of $\frak d_{\xg}\psi_i^k$ in
   $\Lambda^{k+1,N-k}(\xg/\xm)^*$. To prove (ii), we thus have to show that this
   component equals $2(\rho-\mu_i)E^*\wedge\psi_i^k$.

   To do this, we can insert elements of $\xg$ and since both sides have the same
   bidegree we can take $k+1$ elements of $\xp$ and $N-\ell$ elements of $\xk$. Since
   we work modulo $M$ the first $k+1$ elements can either all come from $\xp_+$ or we
   can take $E$ together with $k$ elements of $\xp_+$. Moreover, the last $N-\ell$
   elements can be chosen as $\tilde Y_j=\theta(Y_j)+Y_j$ for $Y_j\in\xg_{r_j}$ for
   some $r_j>0$.  Now since $\mE_i^k\subset\ker(\square_P)$, we get
   $\frak d_P(\beta_s)=0$ and by definition $E\cdot\beta_s=\mu_i\beta_s$ for any
   $s=1,\dots, L$. Again by definition $E\cdot\vol_M=2\rho\vol_M$ which shows that
   $E\cdot\beta^\#_s=(2\rho-\mu_i)\beta^\#_s$ for each $s=1,\dots,L$.  Hence we
   conclude that each of the summands in the formula for $\psi_i^k$ satisfies the
   assumptions of Lemma \ref{lem:delg} with $\nu=2\rho-\mu_i$. In particular,
   $2(\rho-\mu_i)E^*\wedge\psi_i^k$ can produce a non-zero value after insertions of
   elements $\tilde Y_{k+1},\dots,\tilde Y_{k+\ell}$ only if
   $r_{k+1}+\dots+r_{k+\ell}=\nu$. But since $d_K\varphi_i^k$ is the Poisson kernel
   of $d\circ\Phi_i^k$ it follows that the component of $\frak d_{\xg}\psi_i^k$ in
   $\Lambda^{k+1,N-k}(\xg/\xm)^*$ must have the same property. Hence Lemma
   \ref{lem:delg} shows that both side agree upon insertion of appropriate elements
   of $\xg$, which completes the proof of (ii).

   For part (i), we first observe that, up to a sign, we can compute
   $\delta_K\varphi_i^k$ as $*_Kd_K*_K\varphi_i^k$. Now $*_K$ comes from the
   $K$-invariant inner product on $\xq$ for which $\xa$ is orthogonal to
   $\xq\cap(\xg_-\oplus\xp_+)$. Moreover, up to a fixed factor, this inner product
   comes from $\langle\ ,\ \rangle$ on elements of the form $i_k(\alpha)$. In
   particular, up to a fixed factor, we get  $*_Ki_k(\beta_s)=E^*\wedge
   i_k(\beta^\#_s)$ and hence
   $*_K\psi_i^k=\sum_{s=1}^LE^*\wedge i_k(\beta^\#_s)\wedge\beta^\#_s$. 

   Now $d_K$ satisfies a Leibniz rule and formula \eqref{eq:delg} easily implies that
   $d_KE^*=0$. Now for $\alpha\in\Lambda^{0,k-1}(\xg/\xm)^*$ we get
   $0=\langle\beta_s,\frak d_P\alpha\rangle$ and hence $0=\beta^\#_s\wedge\frak
   d_P\alpha$. But now we know that $\Lambda^{0,N}(\xg/\xm)^*$ is one-dimensional and
   by Kostant's theorem the subspace $\ker(\square_P)$ has to be non-trivial and thus
   coincides with $\Lambda^{0,N}(\xg/\xm)^*$.  This shows that $\frak
   d_P:\Lambda^{0,N-1}(\xg/\xm)^*\to \Lambda^{0,N}(\xg/\xm)^*$ is the zero
   map. Together with the Leibniz rule, the above equation thus implies $0=(\frak
   d_P\beta^\#_s)\wedge\alpha$ and since this holds for any $\alpha$, we get $\frak
   d_P\beta^\#_s=0$. Hence any of the forms $i_k(\beta^\#_s)\wedge\beta^\#_s$ satisfies
   the assumptions of Lemma \ref{lem:delg} so we conclude as in the proof of part
   (ii) that after applying $d_K$ and summing up, we obtain a multiple of $E^*$. But
   from above we know that we have to wedge this with $E^*$ to obtain (a multiple of)
   $d_K *_K\varphi_i^k$. Thus $d_K *_K\varphi_i^k=0$ and the proof is complete.
   \end{proof}

  \section{Harmonic Poisson transforms for $\SU(n+1,1)$}\label{3}
  Up to this point, our results were valid for all cases of real rank
    one or in even greater generality. Once one gets into the underlying geometric
    structures and the explicit form of the BGG complex (including the orders of
    operators in the complex), there seems to be no hope for a uniform
    treatment. Hence from this point on, we only discuss the complex case, so
    $G=\SU(n+1,1)$ for the rest of the paper. 
  
\subsection{Geometric setup}
We realize $G=\SU(n+1,1)$ as the stabilizer in $\GL(n+2, \mc)$ of the sesquilinear
form $(z, w) \mapsto z_0\overline{w_{n+1}} + z_{n+1}\overline{w_0} + \sum_{j=1}^n z_j
\overline{w_j}$ for $z$, $w \in \mc^{n+2}$. The maximal compact subgroup $K$ is
isomorphic to $S(\U(n+1)\times \U(1))$, i.e. all pairs $(B, b) \in \U(n+1) \times
\U(1)$ with $b\det(B)= 1$. The parabolic subgroup $P$ is isomorphic to $\CSU(n+1)
\rtimes (\mc^{n} \oplus i\mr)$ and $M$ is given by $S(\U(n) \times \U(1))$, i.e. all
pairs $(B, b) \in \U(n) \times \U(1)$ so that $b^2\det(B)=1$. The Lie algebra $\xg =
\mathfrak{su}(n+1,1)$ is naturally $|2|$-graded with grading components $\xg_{\pm 1}
\cong \mc^n$ and $\xg_{\pm 2} \cong i\mr$, whereas the component $\xg_0$ splits into
the direct sum of $\xa \cong \mr$ and $\xm \cong i\mr \oplus \mathfrak{su}(n)$. Under
these identifications the $M$-action on $\xg_{\pm 1}$ is given by $(B, b) \cdot X =
b^2BX$ for all $X \in \xg_{\pm 1} \cong \mc^n$, whereas the components $\xg_{\pm 2}$
are $M$-invariant.

\textbf{The space $G/K$.} The symmetric space $G/K$ is the complex hyperbolic space
of (complex) dimension $n+1$, which is a K\"{a}hler manifold. The complex structure
on $G/K$ induces a decomposition of
$\Omega^k(G/K, \mc) = \bigoplus_{p+q=k} \Omega^{p,q}(G/K)$ into $(p,q)$-types and the
exterior derivative splits accordingly as $d = \partial + \opartial$. Next, the
K\"{a}hler metric induces an $L^2$-inner product on the space of differential forms
and we define $\delta$, $\partial^*$ and $\opartial^*$ to be the formal adjoints of
$d$, $\partial$ and $\opartial$, respectively. Consequently, we have the Laplacian
$\Delta = d\delta +\delta d$ and similarly
$\Box = \partial\opartial + \opartial\partial$ and its complex conjugate
$\overline{\Box}$, which are all formally self adjoint with respect to the
$L^2$-inner product and satisfy $\Delta = 2\Box = 2\oBox$,
c.f. \cite{ballmann}*{Corollary 5.26}. Additionally, the wedge product with the
K\"{a}hler form $\omega \in \Omega^{1,1}(G/K)$ yields the Lefschetz map $L$ and we
denote by $L^*$ its adjoint with respect to the induced inner product on
$\Lambda^k T^*(G/K)$ and call it the Co-Lefschetz map. These operators are linked via
the K\"{a}hler relations $[L^*, \partial] = i\opartial^*$ and
$[L^*, \opartial] = - i\partial^*$, c.f. \cite{ballmann}*{Proposition 5.22}. A
differential form $\alpha \in \Omega^k(G/K)$ is said to be \emph{primitive} if
$L^*\alpha = 0$, and primitive forms only exist for $0 \le k \le n+1$.

We have seen that the images of the transforms from Definition \ref{def-our-Poisson}
consist of coclosed and harmonic differential forms. Using the K\"{a}hler relations
we can prove an equivalent description of these two properties which will
significantly simplify necessary computations.

\begin{lemma}\label{lem_kaehler_equivalence_harmonic_partial}
 Let $\alpha \in \Omega^{p,q}(G/K)$ be a primitive differential form. Then $\alpha$
 is coclosed if and only if $d\alpha$ is primitive. Furthermore, if $\alpha$ is
 coclosed then $\alpha$ is an eigenform of $\Delta$ with eigenvalue $\lambda \in \mc$
 if and only if $L^*\opartial\partial\alpha = \frac{\lambda}{2i}\alpha$.
\end{lemma}
\begin{proof} 
 By primitivity of $\alpha$ the K\"{a}hler relations simplify to $\partial^*\alpha =
 iL^*\opartial\alpha$ and $\opartial^*\alpha = -iL^*\partial\alpha$, implying the
 first equivalence. Next, $\Delta\alpha = \lambda \alpha$ if and only if $\Box\alpha
 = \frac{\lambda}{2}\alpha$, which by coclosedness of $\alpha$ is equivalent to
 $\partial^*\partial\alpha = \frac{\lambda}{2}\alpha$. Applying $[L^*, \opartial] = -
 i\partial^*$ shows that this relation is equivalent to $\frac{\lambda}{2i}\alpha =
 L^*\opartial\partial\alpha - \opartial L^*\partial\alpha = L^*\opartial\partial
 \alpha$, where the last equality uses the first part.
\end{proof}

\textbf{The space $K/M$.}  In the next step we describe the geometry of the
homogeneous space $G/P$ in terms of its underlying pseudo-Hermitian structure induced
by the isomorphism $K/M \cong G/P$. The group $M$ acts trivially on the grading
component $\xg_2 \cong i\mr$, so each element induces a $K$-invariant $1$-form
$\vartheta$ on $K/M$ which is a contact form and unique up to multiples. We define
$\zeta \in \mathfrak{X}(K/M)$ to be the Reeb vector field corresponding to
$\vartheta$ which is characterized by $\vartheta(\zeta) = 1$ and
$\iota_{\zeta}d\vartheta = 0$. Then $T(K/M)$ is the direct sum of contact subbundle
$H := \ker(\vartheta)$ and the line bundle spanned by $\zeta$, which is isomorphic to
the quotient $Q := T(K/M)/H$. Since $H$ is associated to the $M$-module
$\xg_{-1} \cong \mc^n$ it carries a canonical $K$-invariant complex structure $J$
which induces a decomposition $H \otimes \mc = H^{1,0} \oplus H^{0,1}$ into
holomorphic and antiholomorphic parts and accordingly a decomposition
$\Lambda^k H^* \otimes \mc = \bigoplus_{p+q=k} \Lambda^{p,q} H^*$ into
$(p,q)$-types. In particular, we obtain a splitting
\begin{align}\label{eq:splitting_exterior_power_boundary}
 \Lambda^k T^*(K/M) \otimes \mc = \bigoplus_{p+q=k} \Lambda^{p,q} H^* \oplus \bigoplus_{p+q=k-1} Q^* \otimes \Lambda^{p,q} H^*
\end{align}
for all $0 \le k \le n$, and since $[\cl_{\zeta}, J] = 0$ this decomposition is
invariant under the Lie derivative with respect to the Reeb vector field.

Since $\vartheta$ is a contact form, the restriction
$\omega_H\in \Gamma(\Lambda^{1,1} H^*)$ of $d\vartheta$ is nondegenerate and thus
$g_H(\ , \ ) := \omega_H(\ , J\ )$ defines a bundle metric on $H$. Denote by
$\langle \ , \ \rangle_H$ the induced pointwise inner product on $\Lambda^{*} H^*$
and by $\langle\!\langle \ , \ \rangle\!\rangle_H$ its induced $L^2$-inner product on
$\Gamma(\Lambda^* H^*)$. The $(n,n)$-form $\vol_H := \frac{1}{n!} \omega_H^n$ is
normed with respect to $\langle \ , \ \rangle_H$ and thus we define the
\emph{$H$-Hodge star operator}
$\ast_H \colon \Lambda^{p,q} H^* \to \Lambda^{n-q,n-p} H^*$ via the relation
$\alpha \wedge\ast_H\beta = \langle \alpha, \overline{\beta}\rangle_H \vol_H$ for all
$\alpha \in \Lambda^{q,p} H^*$.

The inclusion of $H$ into $T(K/M)$ dualizes to the projection $\pi_H \colon \Lambda^k
T^*(K/M) \to \Lambda^k H^*$ and we define $d_H \colon \Gamma(\Lambda^k H^*) \to
\Gamma(\Lambda^{k+1} H^*)$ via $d_H := \pi_H \circ d$. Since both $H^{1,0}$ and
$H^{0,1}$ are involutive subbundles of $T(K/M) \otimes \mc$ we can split the complex
linear extension of $d_H$ into the sum of $\partial_H$ and $\opartial_H$, where the
first and second operator maps sections of $\Lambda^{p,q} H^*$ to sections of
$\Lambda^{p+1,q} H^*$ and $\Lambda^{p,q+1} H^*$, respectively.  Next, define
$\delta_H$ as the (formal) $L^2$-adjoint of $d_H$, which again splits into $\delta_H
= \partial_H^* + \opartial_H^*$ according to $(p,q)$-types. Finally, we obtain the
$H$-Laplacian $\Delta_H := d_H\delta_H + \delta_H d_H$ as well as the operator $\Box_H
:= \partial_H\partial_H^* + \partial_H^*\partial_H$ and its complex conjugate.

On sections of $\Lambda^{p,q} H^*$, define the map $L_H$ as the wedge product with
$\omega_H$ and $L_H^*$ as its adjoint with respect to $\langle \ , \ \rangle_H$. It
is immediate that $L_H$ and $L_H^*$ raises respectively lowers the bidegree by
$(1,1)$, and a chase through the definitions show that $L_H^*$ coindices with the
Levi-bracket under the identification of $Q$ with the line subbundle of $T(K/M)$
spanned by the Reeb vector field. Next, by nondegeneracy of $\omega_H$ we see that
$L_H$ is injective on $\Lambda^{p,q}H^*$ for all $0 \le p + q < n$, which by duality
also implies that $L_H^*$ is surjective on this bundle. Moreover, via the duality
given by the $H$-Hodge star we deduce that $L_H$ is surjective and $L_H^*$ is
injective on $\Lambda^{p,q}H^*$ for $p+q \ge n+1$. In the remaining case $p+q = n$
the kernels of $L_H$ and $L_H^*$ coincide.

\begin{prop}\label{prop_operators_pseudo_hermitian}
\begin{enumerate}[(i)]
 \item The assignment $x \mapsto L_H^*$, $y \mapsto L_H$ and $h \mapsto \sum_{p,q}
   (n-(p+q))P_{p,q}$ defines a representation of $\mathfrak{sl}_2(\mc)$ by vector
   bundle homomorphisms on $\Lambda^{*,*} H^*$, where $P_{p,q}$ is the
   projection from $\Lambda^{p+q} H^*$ onto its $(p,q)$-type.
 \item The Lie derivative along $\zeta$ satisfies $\cl_{\zeta}^* = - \cl_{\zeta}$ and
   commutes with $L_H$, $\partial_H$, $\opartial_H$ as well as with their formal
   adjoints.
 \item The operators $\partial_H$, $\opartial_H$, $\partial_H^*$ and $\partial_H^*$ all square to zero and satisfy
 \begin{align*}
  \partial_H \opartial_H + \opartial_H \partial_H + L_H\cl_{\zeta} &= 0, &
  \partial_H^*\opartial_H^* + \opartial_H^*\partial_H^* - \cl_{\zeta}L_H^* &= 0.
 \end{align*}
Moreover, $L_H$ commutes with $\partial_H$ and its complex conjugate, whereas $L_H^*$ commutes with $\partial_H^*$ and its complex conjugate.
\item We have the ``K\"{a}hler relations''
 \begin{align*}
 [L_H^*, \partial_H] &= i \opartial_H^*, & [L_H^*, \opartial_H] &= -i \partial_H^*, & [L_H, \partial_H^*] &= i \opartial_H, & [L_H, \opartial_H^*] &= -i \partial_H.
\end{align*}
\item The partial Laplacians $\Box_H$ and $\overline{\Box}_H$ are related to $\Delta_H$ on $\Gamma(\Lambda^{p,q}H^*)$ via
\begin{align*}
 2 \Box_H &= \Delta_H - i(n-(p+q))\cl_{\zeta},  & 2 \overline{\Box}_H &= \Delta_H + i(n-(p+q))\cl_{\zeta}.  
\end{align*}
In particular, for $p+q = n$ we have $\Delta_H = 2\Box_H = 2\overline{\Box}_H$. Moreover, all Laplacians preserve primitivity.
\end{enumerate}
\end{prop}
\begin{proof}
(i) The only nontrivial relation is $[x,y] = h$, which can be easily obtained
   using the description of $L_H$ and $L_H^*$ in terms of a local orthonormal frame of
   $H$.
   
(ii) Since $\iota_{\zeta}\omega_H = 0$ we deduce that $[\cl_{\zeta}, L_H] = 0$, and from $g_H(\ , \ ) = \omega_H(\ ,J\ )$ together with $[\cl_{\zeta}, J] = 0$ we obtain $\cl_{\zeta}^* = -\cl_{\zeta}$, which in turn implies $[\cl_{\zeta}, L_H^*] = 0$ by taking adjoints. Finally, from \cite{julg_kasparov}*{Proposition 5.10} we know that $[\cl_{\zeta}, d_H] = 0$, so splitting into $(p,q)$-types and taking adjoints yields the rest.

(iii) The first relations follow from expanding $d^2 = 0$ together with taking adjoints, c.f. \cite{rumin}*{eqn. (5), (6)}. For the second part $d_H\omega_H = 0$ immediately implies $[L_H, \partial_H] = 0$ and similarly for its complex conjugate, and the rest follows by taking adjoints. 

(iv) Follows from \cite{julg_kasparov}*{Lemma 5.11} by decompositon into $(p,q)$-types.

(v) Using the relations from (iv) it is easy to see that
   $\partial_H\opartial_H^* = - \opartial_H^*\partial_H$ and similarly for their
   complex conjugates, from which $\Delta_H = \Box_H + \oBox_H$ follows by inserting
   the definition of $\Delta_H$. On the other hand, \cite{rumin}*{Proposition 1}
   shows that $\Box_H - \oBox_H = i((p+q)-n)\cl_{\zeta}$ implying the formulae for
   $\Box_H$ and $\oBox_H$. Finally, the K\"{a}hler relations in (iii) yield
   $[L_H^*, \Box_H] = i \cl_{\zeta} L_H^*$, implying that $\Box_H$ preserves
   primitivity. Since $L_H^*$ is a real operator the rest follows by complex
   conjugation.
\end{proof}

In the next step we describe the BGG-complex on $G/P$ in terms of the operators
induced by the pseudo-Hermitian structure. The group $P$ acts on $\xg_2$ by scalars,
so the subbundle $H = \ker(\vartheta)$ is $G$-invariant. Furthermore, by naturality
of the exterior derivative the group $G$ acts on the restriction $\omega_H$ of
$d\vartheta$ by multiplication with positive smooth functions. In particular, for all
$0 \le k \le 2n$ kernel of $L_H^*|_{\Lambda^k H^*}$ is an invariant subbundle which
we denote by $\Lambda^k_0 H^*$. The definition of the Kostant codifferential easily
implies that $\frak d^*\alpha$ coincides with $\vartheta \wedge L_H^*\alpha$ up to a
multiple. Therefore, the homology bundle $\ch_k$ is isomorphic to $\Lambda^k_0 H^*$
for $k \le n$ and to a quotient of $\Lambda^{k-1}H^*\otimes Q^*$ for $k \ge n+1$,
c.f.\ \cite{cap_harrach_julg}*{Proposition 3.1}. Additionally, on the fibre of $H$
over the origin of $G/P$ the group $P_+$ acts trivially, wheras $G_0$ acts by complex
linear maps, implying that the complex structure $J$ on $H$ and thus the
decomposition of $\Lambda^k H^*$ into $(p,q)$-types is $G$-invariant. Thus, for $k
\le n$ we obtain an induced decomposition of $\ch_k \otimes \mc = \bigoplus_{p+q=k}
\ch_{p,q}$ into $(p,q)$-types, and it turns out that each of the $G$-invariant
bundles $\ch_{p,q}$ is irreducible. Accordingly, for $1 \le k \le n$ the $k$-th
BGG-operator $D = D_k \colon \Gamma(\ch_{k-1}) \to \Gamma(\ch_{k})$ splits into the
sum $D = \cd + \overline{\cd}$, where $\cd$ and $\ocd$ maps sections of $\ch_{p,q}$
to sections of $\ch_{p+1,q}$ and $\ch_{p,q+1}$, respectively.

\begin{thm}\label{thm_BGG_operator}
 Let $0\le p, q \le n$ and define $\lambda_{p,q} := n+1-(p+q)$. 
 \begin{enumerate}[(i)]
\item For $p+q <
n$ the partial BGG-operator is given by $\cd = \partial_H -i \lambda_{p,q}^{-1} L_H\opartial_H^*$.
\item The BGG-operator $D_{n+1} \colon \Gamma(\ch_n) \to \Gamma(\ch_{n+1})$ is given by
\begin{align*}
 D_{n+1} = \vartheta \wedge\left(\cl_{\zeta} + i d_H\opartial_H^* - id_H\partial_H^*\right).
\end{align*}
In particular, the image of the operator $\iota_{\zeta}D_{n+1}$ consists of primitive $n$-forms.
 \end{enumerate}
\end{thm}
\begin{proof}
\cite{rumin}*{Proposition 4}.
\end{proof}

Via the identification of $\ch_{p,q}$ as a subbundle of $\Lambda^{p+q}_0 H^* \otimes
\mc$, we obtain a pointwise Hermitian product on the holonomy bundles well as an
induced $L^2$-Hermitian product, which are both $K$-invariant. Denote by $D^* = D^*_k
\colon \Gamma(\ch_k) \to \Gamma(\ch_{k-1})$ the $L^2$-adjoint of the $k$-th
BGG-operator. For $k \le n-1$ we can split $D^*$ into $\cd^* + \overline{\cd}^*$,
where the first and second operator lowers the first and second degree,
respectively. Using the formula for the BGG-operator in Theorem
\ref{thm_BGG_operator} we immediately deduce that $\cd^*$ corresponds to
$\partial_H^*$ under the identification $\ch_{p,q} \cong \Lambda^{p,q}_0 H^*$, and
similarly for their complex conjugates. In particular, Proposition
\ref{prop_operators_pseudo_hermitian}(iii) implies that $\cd^*$ and $\ocd^*$
anticommute.

\begin{cor}\label{cor_commutator_relation_BGG_adjoint}
 For $k \le n-1$ the BGG-operators and their adjoints satisfy on
 $\Gamma(\Lambda^{k}_0 H^*)$ the relations
 \begin{align*}
 \lambda_k\ocd^*_k \cd_k + (\lambda_k -1) \cd_{k-1}\ocd^*_{k-1} &= 0, & \lambda_k
 \cd^*_{k}\ocd_{k} + (\lambda_k -1) \ocd_{k-1} \cd^*_{k-1} &= 0,
 \end{align*}
 where $\lambda_k = n+1-k$.
\end{cor}
\begin{proof}
 Express $\cd$ and $\ocd$ using Theorem \ref{thm_BGG_operator} and apply Proposition
 \ref{prop_operators_pseudo_hermitian}.
\end{proof}

\subsection{The distinguished harmonic Poisson transform}
The invariant structures on the tangent bundles of $G/K$ and $G/P$ induce
corresponding structures on the subbundles $T_K$ and $T_P$ of $T(G/M)$, which yield
further properties of the $G$-equivariant operators from Definition
\ref{def-our-Poisson}.

The bunde $T_P = \ker(\pi_K)$ is assciated to the same $M$-representation as
$T(K/M)$, implying that each tensorial object on the latter space induces analogous
objects on $T_P$, which we denote by the same symbol. Explicitly, $T_P$ contains a
rank $2n$-bundle $H$, which is endowed with an invariant complex struture $J$ as well
as a bundle metric. The fibre $H_{eM}$ of $H$ at the origin $eM \in G/M$ is
isomorphic to $\mc^n$, and the $M$-invariant inner product on $H_{eM}$ is a multiple
of the standard Hermitian inner product. Moreover, there is an invariant section
$\omega_H \in \Gamma(\Lambda^2 T_P^*)$, and wedging with $\omega_H$ defines a
tensorial map $L_H$ with adjoint $L_H^*$. The kernel $\Lambda^k_0 H^*$ of the
restriction of $L_H^*$ to $\Lambda^k H^*$ is an $G$-invariant subbundle, whose
complexification decomposes into $(p,q)$-types. From the description of the homology
bundles on $G/P$ it follows that the subbundle $\ker(\Box_P) \subset \Lambda^k T^*_P$
coindices with $\Lambda^k_0 H^*$ for $0 \le k \le n$.

On the other hand, we have seen in Section \ref{sec:our-Poisson} that the grading
element $E \in \xa$ induces an invariant vector field
$\zeta_E \in \Gamma(T_K)\subset \mathfrak{X}(G/M)$ as well as an invariant $1$-form
$E^* \in \Omega^1(G/M)$ of bidegree $(1,0)$. Denoting by $E^*_{\mc}:T_K\to\mc$ the complex linear extension of $E^*:T_K\to\mr$,
$\ker(E^*_\mc) \subseteq T_K$ is a $G$-invariant subbundle whose underlying
representation is isomorphic to $\mc^n$. Furthermore, the pullback of the metric and
the complex structure on $G/K$ along the projection $\pi_K \colon G/M \to G/K$ yields
an invariant bundle metric as well as an invariant complex structure on $T_K$,
respectively, and the restriction of the metric to $\ker(E^*_\mc)$ corresponds to the
standard Hermitian inner product on $\mc^n$. A vector field $\xi \in \Gamma(T_K)$
determines sections $\xi^{1,0}$ and $\xi^{0,1}$ of the holomorphic and the
anti-holomorphic part of $T_K\otimes\mc$.  The Co-Lefschetz map on $G/K$ is tensorial
and induces a $G$-invariant complex linear map $L_K^*$ on
$\Lambda^* T^*_K\otimes\mc$. Explicitly, if $\{\xi_{\ell}\}$ is a complex orthonormal
frame of $\ker(E^*_\mc)\subset T_K$, then we can write $L_K^*$ as a linear
combination of $\iota_{\zeta_E^{1,0}}\iota_{\zeta_E^{0,1}}$ and
$\sum_{\ell}\iota_{\xi_{\ell}^{1,0}}\iota_{\xi_{\ell}^{0,1}}$. Observe that after
complexification $E^*_\mc$ induces two sections $E^{*(1,0)}$ and $E^{*(0,1)}$ of
$(T_K\otimes\mc)^*$ which vanish upon insertion of anti-holomorphic respectively
holomorphic vector fields.

For $0 \le p, q \le n$ with $0 \le p+q=k \le n$ we denote the Poisson transform as in
Definition \ref{def-our-Poisson} corresponding to the irreducible subbundle
$\ch_{p,q}$ of $\ch_k \otimes \mc$ by $\Phi_{p,q} \colon \Omega^k(G/P, \mc) \to
\Omega^k(G/K, \mc)$ with induced operator $\underline{\Phi}_{p,q} \colon
\Gamma(\ch_{p,q}) \to \Omega^k(G/K, \mc)$. From the construction of the underlying
map $\phi_{p,q} \in \Hom_G(\Lambda^{p,q}_0 H^*, \Lambda^{p+q} T_K^* \otimes \mc)$ it
is obvious that $\phi_{p,q}$ preserves the $(p,q)$-type and is conformal with respect
to the bundle metrics on $H$ and $\ker(E^*_\mc)$.

\begin{prop}\label{prop_Poisson_complex_case_p_q}
The maps $\underline{\Phi}_{p,q} \colon \Gamma(\ch_{p,q}) \to \Omega^{p+q}(G/K, \mc)$ satisfy the following:
\begin{enumerate}[(i)]
\item The operator $\underline{\Phi}_{p,q}$ has values in $\Omega^{p,q}(G/K)$ and satisfies $\overline{\underline{\Phi}_{p,q}(\sigma)} = \underline{\Phi}_{q,p}(\overline{\sigma})$.
\item The image of $\underline{\Phi}_{p,q}$ consists of primitive, coclosed and harmonic differential forms.
\item There exist complex numbers $c_{p,q} \in \mc$ with $c_{p,q} \neq 0$ so that \begin{align*}
c_{p,q} \cdot \partial \circ \underline{\Phi}_{p-1,q} 
&= (n+1-(p+q)) \underline{\Phi}_{p,q} \circ \cd, \\
\overline{c_{q,p}} \cdot \opartial \circ \underline{\Phi}_{p,q-1} 
&= (n+1-(p+q)) \underline{\Phi}_{p,q} \circ \ocd.
\end{align*}
Moreover, these complex numbers satisfy $c_{p,q}\overline{c_{q,p-1}} = c_{p,q-1} \overline{c_{q,p}}$. 
\end{enumerate}

\end{prop}
\begin{proof}
(i) The operator $\phi_{p,q}$ preserves $(p,q)$-types and satisfies
   $\overline{\phi_{p,q}} = \phi_{q,p}$ by construction, thus the claim follows
   Theorem \ref{thm_uniqueness}(ii).
 
 (ii) In view of Theorems \ref{thm_uniqueness} and \ref{thm_PT_rank_1} it remains to
   prove that $L^* \circ \underline{\Phi}_{p,q} = 0$, which by tensoriality is
   equivalent to $L_K^*\circ \phi_{p,q} = 0$. But by construction the vector fiels
   $\zeta_{E}^{1,0}$ and $\zeta_E^{0,1}$ insert trivially into $\phi_{p,q}(\alpha)$
   for all $\alpha \in \Gamma(\Lambda^{p,q}_0 H^*)$, whereas
   $\sum_{\ell} \iota_{\xi_{\ell}^{1,0}}\iota_{\xi_{\ell}^{0,1}}\phi_{p,q}(\alpha)$
   is a multiple of $\phi_{p-1,q-1}(L_H^* \alpha) = 0$ for each complex orthonormal
   frame $\{\xi_{\ell}\}$ of $\ker(E^*_\mc)\subset T_K$.
 
 (iii) By definition, the bundles $\Lambda^{p,q} T_K^*$ are induced by
   the restriction to $M$ of a representation of $K$ and the bundle map $L_K^*$ is
   induced by a $K$-equivariant map. Hence also $\Lambda^{p,q}_0 T_K^*$ is induced by
   the restriction to $M$ of a representation of $K$, which is well known to be
   $K$-irreducible. 
 
 From (i) and (ii), we see that that $\uPhi_{p,q} \circ \cd$ is an equivariant
 operator $\Gamma(\ch_{p-1,q}) \to \Omega^{p,q}(G/K)$ whose image consists of
 primitive differential forms and hence lies in $\Gamma(\Lambda^{p,q}_0
   T_K^*)$. Moreover, from Proposition \ref{prop_poisson_on_irreducible} we know
 that this operator is again induced by a Poisson transform and hence comes from an
 $M$-equivariant map $\psi \in \Hom_M(\Lambda^{p-1,q}_0 H^*_{eM}, \Lambda^{p,q}_0
 (T_K^*)_{eM})$. So $\psi$ is defined on an irreducible representation of $M$ and has
 values in an irreducible representation of $K$, and from \cite{baldoni}*{Theorem
   4.4} we know that the dimension of the space of homomorphisms is at most one under
 these assumptions. Now $E^{*(1,0)}\wedge \phi_{p-1,q}$ is an
   $M$-equivariant map between the same spaces as $\psi$ and thus must be a
   multiple. In particular, if $\varphi_{p,q}$ is the Poisson kernel underlying
 $\uPhi_{p,q}$, then this shows that the Poisson kernel underlying $\uPhi_{p,q} \circ
 \cd$ is given by $c_{p,q}E^{*(1,0)}\wedge \varphi_{p-1,q}$ for $c_{p,q}
 \in \mc$, and it is easy to see that $c_{p,q} \neq 0$.
 
  Additionally, we know from Theorem \ref{thm_PT_rank_1}(ii) applied for
    $\mk=\mc$ that the Poisson kernel underlying $d \circ \underline{\Phi}_{p,q}$ is
  given by $2(\rho - \mu_{p,q}) E^{*(1,0)}\wedge \varphi_{p,q}$, where
  $2\rho$ and $\mu_{p,q}$ are the $\ad(E)$-eigenvalues of $\Lambda^{2n+1} \xp_+$ and
  $\Lambda^{p,q}_0 (H^*)_{eM}$, respectively. Explicitly, we have $\rho = n+1$ and
  $\mu_{p,q} = p+q$, and decomposing the above kernel into types implies that
  $\partial \circ \underline{\Phi}_{p,q}$ is induced by the kernel $2(n+1-(p+q))
  E^{*(1,0)}\wedge \varphi_{p,q}$. Comparing this to the kernel of
  $\underline{\Phi}_{p-1,q} \circ \cd$ yields the first relation, and the second
  follows from (i) by complex conjugation. Finally, the relation on the coefficients
  follows from anticommutativity of $\partial$ and $\opartial$ and of the partial
  BGG-operators.
\end{proof}

\begin{rem}
  In \cite{cap_harrach_julg}*{Theorem 4.6} and \cite{harrach_thesis} a family of
  harmonic Poisson kernels was constructed whose induced $G$-equivariant maps
  $\underline{\Phi}_{p,q} \colon \Gamma(\ch_{p,q})\to\Omega^{p,q}(G/K, \mc)$ satisfy
  the properties of Proposition \ref{prop_Poisson_complex_case_p_q} with
  $c_{p,q} = 2i(n+2-p)$. (The last part follows by comparing equations (4.18) and
  (4.19) in the proof of \cite{cap_harrach_julg}*{Theorem 4.10}). However, the
  construction of these kernels needed a significant amount of explicit computations,
  which we avoided completely using our general approach.
\end{rem}

Now we can construct a family of real operators $\uPhi_k \colon \Gamma(\ch_k) \to
\Omega^k(G/K)$ for $0 \le k \le n$ via suitable linear combinations of the operators
$\uPhi_{p,q}$ in such a way that we obtain a chain map between the BGG-complex on
$G/P$ and the deRham complex on $G/K$.

\begin{thm}\label{thm_Poisson_complex_space}
For all $0 \le p, q \le n$ with $0 \le p+q \le n$ define $\lambda_{p,q} \in \mc$ via
$\lambda_{0,0} := \frac{1}{\vol(K/M)}$ and
\begin{align*}
 \lambda_{p,q} := \lambda_{0,0}\frac{n!}{(n+1-(p+q))!} \prod_{i=1}^p c_{i,q}^{-1} \prod_{j=1}^q \overline{c_{j,0}}^{-1},
\end{align*}
where $c_{p,q} \neq 0$ are the constants appearing in Proposition \ref{prop_Poisson_complex_case_p_q}(iii). Then for all $0 \le k \le n$ the image of the operators $\uPhi_k \colon \Gamma(\ch_k) \to \Omega^k(G/K)$ defined by $\uPhi_k := \sum_{p+q=k} \lambda_{p,q} \uPhi_{p,q}$
consist of primitive, coclosed and harmonic differential forms. Moreover, for all $1 \le k \le n$ they satisfy the relation $d \circ \uPhi_{k-1} = \uPhi_k \circ D_{k}$.
\end{thm}
\begin{proof}
 
Using the relations $c_{p,q}\overline{c_{q,p-1}} = c_{p,q-1} \overline{c_{q,p}}$ it is easy to see that the coefficents satisfy
\begin{align}\label{eq:relation_lambda}
 (n+1-(p+q))\lambda_{p,q} = \lambda_{p+1,q}c_{p+1,q} = \lambda_{p,q+1} \overline{c_{q+1,p}},
\end{align}
which by part (iii) of the previous Proposition imply that these operators are chain maps.
\end{proof}

The above construction of real Poisson transforms provides a geometric
  connection between the differential calculi on $G/P$ and $G/K$ in each degree. For
  the following developments, the cases of main interest are the middle degrees $n$
  and $n+1$. In the these degrees, we have the diagram
\begin{center}
\begin{tikzcd}
   \dotsc \arrow{r}{D_{n-1}} & \Gamma(\ch_{n-1}) \arrow{r}{D_n}
   \arrow{d}{\underline{\Phi}_{n-1}} & \Gamma(\ch_n) \arrow{r}{D_{n+1}} \arrow[rd,
     dashed]{}{d\underline{\Phi}_n} \arrow{d}{\underline{\Phi}_{n}} &
   \Gamma(\ch_{n+1}) \arrow{r} \arrow{d}{\underline{\Phi}_{n+1}}& \dotsc \\ \dotsc
   \arrow{r} & \Omega_{\Delta}^{n-1}(G/K) \arrow{r}{d} & \Omega_{\Delta}^n(G/K)
   \arrow{r}{d} & \Omega_{\Delta}^{n+1}(G/K) \arrow{r} & \dotsc
\end{tikzcd}
\end{center}
where we denote by $\Omega^k_{\Delta}(G/K)$ the space of harmonic $k$-forms
on $G/K$. Now it is well known that the space $L^2\Omega^{n+1}_{\Delta}(G/K)$ of
$L^2$-harmonic $(n+1)$-forms on $G/K$ is a realization of the direct sum of all
discrete series representations of $G$ with trivial infinitesimal character (see
\cite{pedon_complex} and Section \ref{sec_DSR}). Our aim for the rest of the paper is
to show that the Poisson transform $d\underline{\Phi}_n$ maps into the space of
$L^2$-harmonic forms and has dense image there. On the one hand, this leads to a
realization of (a dense subspace in) the sum of discrete series representations with
trivial infinitesimal character on a quotient of the space
$\Gamma(\ch_n)$ on the boundary sphere. On the other hand, as we will prove in
Theorem \ref{prop_L_2_general} that we obtain an induced exact sequence
  \begin{center}
\begin{tikzcd}
   \dotsc \arrow{r}{D_{n-1}}  & \Gamma(\ch_{n-1}) \arrow{r}{D_n} & \Gamma(\ch_n) \arrow{r}{d\underline{\Phi}_n} \arrow{r} &  L^2\Omega_{\Delta}^{n+1}(G/K) 
\end{tikzcd}
\end{center}
for which the last operator has dense image. So one can ``cut off'' (a completion of)
the BGG complex in the middle degree by the sum of discrete series representations
with trivial infinitesimal character which is important for applications towards the
Baum-Connes conjecture.

\section{Boundary asymptotics of Poisson transforms.}
  It will be of crucial importance for us to prove that in some cases our
  Poisson transforms have $L^2$-values. Thus we have to understand the behavior of
  these forms towards infinity and in order to do this, we use the fact that $G/P$
  can be attached to $G/K$ as a ``boundary at infinity''. This is based on the
  Poincar\'{e} ball model of complex hyperbolic space.

Recall that we realized the group $G$ as the stabilizer of a Lorentzian Hermitian
form on $ \mc^{n+2}$. The maximal compact subgroup $K$ corresponds to the stabilizer
of a negative line $\ell_-$, whereas the parabolic subgroup $P$ is the stabilizer of
a null line, which shows that $G/K$ and $G/P$ are isomorphic to the sets of negative
and null lines in $\mc^{n+2}$, respectively. To see how to realize $G/P$
as the boundary at infinity of $G/K$ fix a vector $v_- \in \ell_-$ of length $-1$ and
consider the map which sends $v \in \ell_-^{\perp}$ to the complex line in
$\mc^{n+2}$ generated by $v + v_-$. This restricts to a diffeomorphism from the open
unit ball in $\ell^{\perp}_-$ to the space of negative lines, whereas the unit sphere
in $\ell^{\perp}$ is mapped diffeomorphically onto the space of null lines. If we
identify $\ell_-^{\perp}$ with $\mc^{n+1}$ this shows that $G/K$ is diffeomorphic to
the open unit ball $B \subset \mc^{n+1}$, and the invariant Hermitian metric
corresponds to the Bergman metric
\begin{align} \label{eq:Bergman}
 h(z)(\xi, \eta) = \frac{\langle \xi, \eta\rangle}{1-|z|^2} + \frac{\langle \xi,
   z\rangle \langle z, \eta \rangle}{(1 -|z|^2)^2}
\end{align}
for all $z \in B$ and $\xi$, $\eta \in T_zB$, where $\langle \ , \ \rangle$ denotes
the standard Hermitian metric on $\mc^{n+1}$. In this picture the parabolic geometry
$G/P$ can be naturally viewed as the boundary sphere $\partial B \cong S^{2n+1}$
with the CR structure inherited as a hypersurface. This is a model case
  for compactifications of K\"ahler-Einstein metrics, see \cite{biquard}.

\subsection{Polar decomposition of complex hyperbolic space}
In the next step we relate differential operators on $K/M$ and $G/K$ in the
Poincar\'{e} ball model via the polar decomposition $\Psi = (r,\chi) \colon B
\setminus \{ 0\} \to (0,1) \times K/M$ given by $r(z) := |z|$ and $\chi(z) :=
r(z)^{-1}z$. For fixed $r \in (0,1)$ the map $\Psi^{-1}_r := \Psi^{-1}(r,-) \colon
K/M \to B \setminus \{0\}$ is a $K$-equivariant isomorphism, implying that the image
$S_r := \Psi^{-1}_r(K/M)$ is a $K$-orbit of $B$ which coincides with the sphere of
radius $r$. Define the distribution $H_r \subset TS_r$ as the image of the
CR-distribution $H \subset T(K/M)$ under $T\Psi^{-1}_r$.

The (complex) standard basis of $\mc^{n+1}$ induces complex coordinates $z = (z_0, \dotsc, z_n)$ of $z \in \mc^{n+1}$. Then $\{z_j, \bz_j\}$ are real coordinates on $B$ and we define the partial derivatives $\partial_{z_j}$ and $\partial_{\bz_j}$ with respect to these coordinates. Define the holomorphic vector field $T := \sum_{j=0}^n z_j \partial_{z_j}$ as well as the $(1,0)$-form $\tau := \sum_{j=0}^n \bz_j dz_j$. We obtain the relations $\tau(T) = \otau(\overline{T}) = r^2$ and therefore $\partial_r = \frac{1}{r}(T + \overline{T})$ and $dr = \frac{1}{2r}(\tau + \otau)$. Furthermore, we fix the pseudo-Hermitian structure $\vartheta \in \Omega^1(K/M)$ so that its corresponding Reeb vector field $\zeta \in \mathfrak{X}(K/M)$ is given by $\zeta(u) = i(T(u) - \overline{T}(u))$ for all $u \in S^{2n+1}$. Using the definition of $\Psi^{-1}_r$ a direct computation shows that $\zeta$ is pushed forward to the vector field $J\partial_r = \frac{i}{r}(T-\overline{T}) \in \mathfrak{X}(B\setminus\{0\})$. 

It is immediate that for $z \in S_r$ the tangent space $T_zB\otimes\mc$ is the direct sum of $H_r\otimes \mc$ and the complex subspace spanned by $T(z)$ and $\overline{T}(z)$. Using this and the description \eqref{eq:Bergman} for the 
Bergman metric on $B$ a direct computation yields that the K\"{a}hler metric $g$ and the K\"{a}hler form $\omega$ are given in terms of the polar decomposition by
 \begin{align*}
  g &= \frac{1}{(1-r^2)^2r^2} \frac{1}{2}\tau \odot \overline{\tau} + \frac{r^2}{1-r^2} \frac{1}{2} \chi^*g_H, &
  \omega &= \frac{1}{(1-r^2)^2r^2}\frac{i}{2} \tau \wedge \overline{\tau} + \frac{r^2}{1-r^2} \frac{1}{2}\chi^*\omega_H.
 \end{align*}
Furthermore, since the distribution $H_r$ was defined as the pullback of $H \subset T(G/P)$ along the angular map every $(p,q)$-form $\alpha$ on $B\setminus \{0\}$ can be decomposed as
\begin{align*}
 \alpha = \chi^*\sigma_{p,q} + \tau \wedge \chi^*\sigma_{p-1,q} + \overline{\tau} \wedge \chi^*\sigma_{p,q-1} + \tau \wedge \overline{\tau} \wedge \chi^*\sigma_{p-1,q-1},
\end{align*}
where $\sigma_{s,t}$ is a smooth map $(0,1) \to \Gamma(\Lambda^{s,t} H^*)$. In the
sequel we determine an explicit formula for the image of the operator
$\underline{\Phi}_{p,q}$ in terms of this decomposition, where we will use that this
image consists of primitive, coclosed and harmonic differential forms. In view of
Lemma \ref{lem_kaehler_equivalence_harmonic_partial} we need to determine formulae
for $\partial$, $\opartial$ and the Co-Lefschetz map $L^*$.

\begin{prop}\label{prop_polar_differential}
 \begin{enumerate}[(i)]
 \item The $1$-form $\tau$ is antiholomorphic and satisfies
   $\opartial\tau = ir^2 \chi^*\omega_H - \frac{1}{r^2}\tau \wedge\overline{\tau}$.
  \item The Lie derivative along $T$ satisfies $\cl_T\chi^*\sigma = \frac{1}{2i}\chi^*\cl_{\zeta}\sigma$ 
for all $\sigma \in \Gamma(\Lambda^{*}H^*)$, where $\zeta$ is the Reeb vector field on $S^{2n+1}$. In particular,
$\partial \chi^*\sigma = \chi^*\partial_H\sigma + \frac{1}{2ir^2}\tau \wedge \chi^*\cl_{\zeta}\sigma$. 
  \item If $\sigma_{s,t} \in \Gamma(\Lambda^{s,t} H^*)$, then the Co-Lefschetz map is determined on $\Omega^{p,q}(B)$ by
  \begin{align*}
   L^\ast\chi^*\sigma_{p,q} &= 2r^{-2}(1-r^2) \chi^*L_H^*\sigma_{p,q}\\
   L^\ast(\tau \wedge \chi^*\sigma_{p-1,q}) &= 2r^{-2} (1-r^2) \tau \wedge \chi^*L_H^*\sigma_{p-1,q}\\
   L^\ast(\overline{\tau} \wedge \chi^*\sigma_{p,q-1}) &= 2r^{-2} (1-r^2) \overline{\tau} \wedge \chi^*L_H^*\sigma_{p,q-1}\\
   L^\ast(\tau \wedge \overline{\tau} \wedge \chi^*\sigma_{p-1,q-1}) &= 2r^{-2}(1-r^2) \tau \wedge \overline{\tau} \wedge \chi^*L_H^*\sigma_{p-1,q-1} - 2i r^2(1-r^2)^2 \chi^*\sigma_{p-1,q-1}.
  \end{align*}
 \end{enumerate}
\end{prop}
\begin{proof}
 (i) Follows from a direct computation using the definition of $\tau$.
 
 (ii) Since $\iota_{\partial_r}d\chi^*\sigma = 0$ by naturality of the exterior derivative it follows that $\iota_{T}\chi^*\sigma = -\iota_{\overline{T}}\chi^*\sigma$. On the other hand, the vector field $i(T - \overline{T})$ is $\chi$-related to the Reeb vector field $\zeta$, implying the first formula. For the second formula we use naturality of the exterior derivative to write
\begin{align}\label{eq:partial_pullback}
 d\chi^*\sigma = \chi^*\partial_H\sigma + \chi^*\opartial_H\sigma + \chi^*\vartheta \wedge \chi^*\cl_{\zeta}\sigma.
\end{align}
A direct computation shows that $\chi^*\vartheta = \frac{1}{2ir^2}(\tau - \overline{\tau})$ and after inserting this into \eqref{eq:partial_pullback} and using (i) we decompose the resulting equation according to $(p,q)$-types.

(ii) Using the expression for the K\"{a}hler form yields explicit formulae for the Lefschetz map and thus for $L^*$ by adjointness.
\end{proof}

\subsection{Boundary values of Poisson transforms}
As a first step towards understanding the boundary asymptotics of the values of
$\underline{\Phi}_{p,q}$, we will consider the case of zeros of the initial section
in the case $p+q\leq n$. Given $x \in \partial B$ and a section $\alpha$ of the
homology bundle with $\alpha(x) = 0$, we prove that the differential form
$\underline{\Phi}_{p,q}(\alpha)$ also tends to $0$ along every smooth curve on
$\overline{B}$ that ends in $x$. This will be done via a norm estimate argument, for
which we need to analyse the behaviour of the $G$-equivariant map underlying
$\Phi_{p,q}$ along the fibres of the projection $\pi_K \colon G/M \to G/K$.

In order to do so, reall that for every $w \in \mr$ there is a real line bundle
$\ce(w) = G \times_P \mr(w) \to G/P$ whose sections are densities of weight $w$,
normalized so that $\ce(-2n-1)$ is isomorphic to the bundle of volume densities on
the oriented manifold $G/P$. Given such a density $\sigma \in \Gamma(\ce(w))$ with
corresponding $P$-equivariant map $f \in C^{\infty}(G, \mr(w))^P$ the pullback
$\pi_P^*\sigma \in C^{\infty}(G/M)$ is the smooth function which is induced by
viewing $f$ as an $M$-equivariant map.

\begin{lemma}\label{lem_metric_volume_form_G_M_comparison}
 Let $G = KAN$ be the Iwasawa decomposition of $G$ and for $g \in G$ write $g = \kappa(g) e^{\tau(g)E}\eta(g)$ with $\kappa(g) \in K$, $\eta(g) \in N$ and $\tau(g) \in \mr$, where $E \in \xa$ is the grading element. Put ${\bf 0} := eK$ and let $z = gK \in G/K$ and $u = hP \in G/P$ with $h \in K$. Define the smooth map $\cp \in C^{\infty}(G/M)$ via $\cp(z,u) := e^{-2\tau(g^{-1}h)}$ using the identification of $G/M$ with $G/K \times G/P$. 
 \begin{enumerate}[(i)]
 \item For all $\alpha$, $\beta \in \Gamma(\Lambda^k H^*)$ we have
 \begin{align*}
  \langle \pi_P^*\alpha, \pi_P^*\beta \rangle_M(x,u) = \cp(z,u)^{-k} \langle \pi_P^*\alpha, \pi_P^*\beta \rangle_M({\bf 0},u).
 \end{align*}
In particular, $\vol_M(x,u) = \cp(z,u)^{n+1}\vol_M({\bf 0}, u)$.
\item The pullback $\pi_P^*\sigma$ of a $w$-density $\sigma \in \Gamma(\ce(w))$ satisfies $\pi_P^*\sigma(x,u) = \cp(z,u)^{\frac{w}{2}}\pi_P^*\sigma({\bf 0}, u)$.
 \end{enumerate}
\end{lemma}
\begin{proof}
The element $h_g := h\eta(g^{-1}h)^{-1}e^{-\tau(g^{-1}h)E}h^{-1} \in G$ maps the point $({\bf 0},u)$ to $(x, u)$ and acts on $H_u$ by multiplication with $e^{\tau(g^{-1}h)}$ and on a vector $\eta \in T_u(G/P)$ transversal to $H_u$ by multiplication with $e^{2\tau(g^{-1}h)}$ plus addition of an element in $H_u$. Applying $G$-invariance of $g_M$ and $\vol_M$ as well as $P$-equivariance of the map underlying $\pi_P^*\sigma$ the claims follow immediately.
\end{proof}

In the next step we need to analyse the asymptotics of powers of $\cp$ as these
govern the change of the differential form $|\pi_P^*\alpha|_M \vol_M$ on $G/M$ along
the fibres over $G/P=S^{2n+1}$ for all $\alpha \in \Gamma(\Lambda^k H^*)$. After
expressing the induced $G$-actions on the Poincar\'{e} ball model and its boundary
sphere an explicit computation yields for $z = gK \in B$ and $u = hP \in G/P$ with
$h \in K$ that
\begin{align*}
 \cp(z,u) = \frac{1-|z|^2}{|1-\langle z, u \rangle|^2}.
\end{align*}
In particular, by Lemma \ref{lem_metric_volume_form_G_M_comparison} the change of the
volume form $\vol_M$ is governed by $\cp(z,u)^{n+1}$, which is (a multiple of) the
Poisson-Bergman kernel on the unit ball in $\mc^{n+1}$, c.f. \cite{krantz}*{Section
  3.2.2}. It turns out that the asymptotic behaviour of $\cp$ can be described via
Gaussian hypergeometric functions with specific parameter, and since these will
appear again later we collect important information in the next Proposition.

\begin{prop}\label{prop_hypergeometric}
Let $a$, $b$, $c$ with $c \not \in \mz_{\le 0}$ and for $z \in \mc$ with $|z| < 1$ let
\begin{align*}
 {}_2F_1(a,b,c;z) := \sum_{n=0}^{\infty} \frac{(a)_n(b)_n}{(c)_n} \frac{z^n}{n!}
\end{align*}
be the (Gaussian) hypergeometric function, where $(q)_n$ is the rising Pochhammer symbol. Then
\begin{enumerate}[(i)]
 \item The function ${}_2F_1(a,b,c;z)$ is the unique solution up to multiples of the hypergeometric differential equation
 \begin{align*}
  z(1-z) F''(z) + (c - (a+b+1)z) F'(z) - ab F(z) &= 0
 \end{align*}
which is bounded for $z = 0$.
\item $\frac{d}{dz} {}_2F_1(a,b,c;z) = \frac{ab}{c} {}_2F_1(a+1,b+1,c+1;z)$
 \item If $\re(c-a-b) > 0$ then the limit of ${}_2F_1(a,b,c;z)$ for $z \to 1$ exists and is given by
 \begin{align*}
  {}_2F_1(a,b,c;1) = \frac{\Gamma(c)\Gamma(c-a-b)}{\Gamma(c-a)\Gamma(c-b)}.
 \end{align*}
\item For $\re(c-a-b) = 0$ we have ${}_2F_1(a,b,c;z) \in L^1([0,1])$.
\end{enumerate}
\end{prop}
\begin{proof}
\cite{daalhuis}*{ch. 15.10, 15.5.1, 15.4.20, 15.4.21}
\end{proof}

Using this information on hypergeometric functions we are ready prove the asymptotic behaviour of specific powers of the function $\cp$. 

\begin{prop}\label{prop_integral_poissonkernel_bounded}
For all $0 \le k \le n+1$ we have
\begin{align*}
 I_k(z) := \int_{S^{2n+1}} \cp(z,u)^{n+1-\frac{k}{2}} du = \vol(S^{2n+1}) \cdot {}_2F_1\left(\frac{k}{2}, \frac{k}{2}, n+1;|z|^2\right).
\end{align*}
In particular, $I_k(z)$ is globally bounded for $0 \le k \le n$.
\end{prop}
\begin{proof}
  First, we show that $I_k$ is an eigenfunction of the Laplacian, which can be done via
  a tedious computation of $\Delta \cp(z,u)^{n+1-\frac{k}{2}}$. However, we can avoid
  this by use of representation theory. Explicitly, combining parts (ii) and (iii) of
  Lemma \ref{lem_metric_volume_form_G_M_comparison} the integrals $I_{-k}(z)$
  coincide with a multiple of the Poisson transform $\Phi_0(\sigma_k)$ of the
  $k$-density $\sigma_{k}$, which is defined by the $P$-equivariant function
  $f_{k} \colon G \to \mr$ given by $f_{k}(g\exp(tE)n):= e^{-kt}$ for any
    $g\in K$ and $n\in N$. In particular, using that $\Delta$ is the negative of the
  differential operator $\cc$ induced by the Casimir element, $G$-equivariance of the
  Poisson transform implies that $\Delta I_{k}(z) = -\Phi_0(\cc \cdot
  \sigma_{-k})$. Now $\sigma_{-k}$ is a section of a $1$-dimensional representation
  with highest weight $\lambda := kE^*$ on which $\cc$ acts by the scalar
  $\langle \lambda, \lambda - 2 \rho \rangle$, c.f. \cite{knapp_lie}*{Proposition
    5.28}. Here the inner product $\langle \ , \ \rangle$ on $\xa^*$ is induced by a
  multiple of the Killing form and normed so that the grading element is of unit
  length. Since the lowest form $\rho$ is given by $(n+1)E^*$, we thus obtain that
  $I_k(z)$ is an eigenfunction of the Laplace operator on $G/K$ with eigenvalue
  $k(2(n+1)-k)$.

  The transformation formula for integrals shows that $I_k(z)$ is $K$-invariant and
  thus only depends on the norm of $z$. Define the smooth function
  $f \colon (0,1) \to \mr$ via $f(r(z)) := I_k(z)$, then $f$ is coclosed and
  primitive. Therefore, we can apply Lemma
  \ref{lem_kaehler_equivalence_harmonic_partial}, which shows that $I_k$ being an
  eigenfunction of $\Delta$ is equivalent to
  $L^*\opartial\partial I_k(z) = \frac{k(2(n+1)-k)}{2i}I_k(z)$. Expanding the left
  hand side using the formulae from Proposition \ref{prop_polar_differential} a few
  lines of computation yields that $f$ has to satisfy the ODE
\begin{align}\label{eq:ODE}
 f''(r) + \left(\frac{2n+1}{r} + \frac{2n r}{1-r^2}\right)f'(r) + \frac{k(2(n+1)-k)}{(1-r^2)^2} f(r) = 0.
\end{align}
After a change of variables to $s(r) = r^2$ this can be transformed into a hypergeometric differential equation in $s$ with parameter $a = b = \frac{k}{2}$ and $c = n+1$. Applying \ref{prop_hypergeometric}(i) it follows that the bounded solution of \eqref{eq:ODE} is given by $f(r) = C \cdot {}_2F_1\left(\frac{k}{2}, \frac{k}{2}, n+1;r^2\right)$ with a constant $C \in \mc$. Forming the limit $r \to 0$ yields that $C = \vol(S^{2n+1})$ and thus the claimed formula. Since $c - a - b = n+1-k > 0$ the boundedness claim follows from Proposition \ref{prop_hypergeometric}(iii).
\end{proof}

Using the previous result we can prove a boundary value theorem for the operators $\uPhi_{p,q}$ in the case of trivial initial value. This will also show that the boundary value of a Poisson transform at a boundary point only depends on the value of the initial data in this point.   

\begin{thm}\label{thm_main_zero_boundary_value}
For $0\le p, q \le n$ with $0 \le p+q \le n$ let $\alpha \in \Gamma(\ch_{p,q})$ with
$\alpha(u) = 0$ for $u \in S^{2n+1}$. Let $\gamma \colon [0, a) \to \overline{B}$ be
  a continuous curve with $\gamma(0) = u$ and $\gamma(t) \in B$ for $t > 0$. Then
 \begin{align}
  \lim_{t \to 0} \underline{\Phi}_{p,q}(\alpha)(\gamma(t)) = 0.
 \end{align}
\end{thm}
\begin{proof}
Using the identification of $\ch_{p,q}$ with $\Lambda^{p,q}_0H^*$ we can find a section $\beta$ of the latter bundle so that $\underline{\Phi}_{p,q}(\alpha) = \Phi_{p,q}(\beta)$. Moreover, we have seen that the $G$-equivariant map $\phi_{p,q}$ underlying the Poisson transform $\Phi_{p,q}$ is a conformal map with respect to the inner products on $\Gamma(\Lambda^k_0 H^*)$ and $\Gamma(\Lambda^k T_K^*)$ whose conformal factor $c_{p,q} \in \mr$ is constant on $G/M$. Thus, if we denote by $\| \ \|$ the norm on $\Omega^{p,q}(G/K)$ induced by the K\"{a}hler metric we obtain for all $z \in B$ that
\begin{align}\label{eq:norm_inequality_1}
 \left\|\uPhi_{p,q}(\alpha)(z)\right\| \le c_{p,q} \int_{S^{2n+1}} \|\pi_P^*\beta\|_M(z,u) \ \vol_M(z,u).
\end{align}
Next, we use Lemma \ref{lem_metric_volume_form_G_M_comparison} to move the integrand to the fibre over the origin ${\bf 0} \in B$, which produces the factor $\cp(z,u)^{n+1-\frac{p+q}{2}}$. Furthermore, along the fibre $\pi_K^{-1}({\bf 0})$ the norm of $\pi_P^*\beta$ coincides with the norm of $\alpha$ with respect to the inner product induced by the $K$-invariant Riemannian metric on $S^{2n+1} \cong K/M$. Inserting this into \eqref{eq:norm_inequality_1} yields
\begin{align}\label{eq:norm_inequality_2}
 \left\|\uPhi_{p,q}(\alpha)(z)\right\| \le c_{p,q} \int_{S^{2n+1}} \cp(z,u)^{n+1-\frac{p+q}{2}} \|\alpha(u)\| \ du.
\end{align}
Now let $\varepsilon > 0$ be arbitrary. Since $\alpha(u) = 0$ we can find $\delta > 0$
so that for all $v \in S^{2n+1}$ in the $\delta$-ball $B_{\delta}(u)$ around $u$ we
have $\|\alpha(v)\| < \varepsilon$. Furthermore, since $n+1 - \frac{p+q}{2} > 0$,
$\cp^{n+1-\frac{p+q}{2}}$ tends to $0$ uniformly on every compact subset on
$S^{2n+1}$ away from $u$. Therefore, since the curve $\gamma$ is continuous we can
find $T > 0$ so that for all $t < T$ and all $v \in S^{2n+1}\setminus B_{\delta}(u)$
we have $\cp(\gamma(t),v)^{n+1-\frac{p+q}{2}} < \varepsilon$. Furthermore, we have
seen in Proposition \ref{prop_integral_poissonkernel_bounded} that the integral over
$\cp(z,u)^{n+1-\frac{p+q}{2}}$ is globally bounded by a constant $C \in
\mr$. Therefore, by \eqref{eq:norm_inequality_2} we have for all $t < T$ that
\begin{align*}
 \|\uPhi_{p,q}(\alpha)(\gamma(t))\| &\le c_{p,q} \left(\int_{B_{\delta}(u)} + \int_{S^{2n+1}\setminus B_{\delta}(u)}\right) \cp(\gamma(t),v)^{n+1-\frac{p+q}{2}} \|\alpha(v)\| \ dv \\
 &< \varepsilon c_{p,q}\int_{B_{\delta}(u)} \cp(\gamma(t),v)^{n+1-\frac{p+q}{2}} \ dv + \varepsilon \int_{S^{2n+1}} \|\alpha(v)\| \ dv \\
 &\le \varepsilon c_{p,q}\left(C + \|\alpha\|_{\infty} \vol\left(S^{2n+1}\right)\right),
\end{align*}
where $\| \alpha \|_{\infty}$ denotes the maximum norm of $\alpha$ on the compact
manifold $S^{2n+1}$. Since $\varepsilon$ was arbitrary the claim follows.
\end{proof}

\section{Poisson transform on isotypic components}
Our next aim is to analyse the image of the operators $\underline{\Phi}_{p,q} \colon
\Gamma(\ch_{p,q}) \to \Omega^{p,q}(G/K)$ for all $0 \le p+q \le n$ defined in
Proposition \ref{prop_Poisson_complex_case_p_q}. In order to do so, recall that the
Peter-Weyl theorem (c.f. \cite{knapp}*{Theorem 1.12}) implies that for a
representation of a compact Lie group $K$ on a Hilbert space $\mw$ the subspace
$\mw^{fin}$ of $K$-finite vectors is dense. Let $\hat{K}$ be the set of all
isomorphism classes of finite dimensional irreducible $K$-representations and for
$\mv \in \hat{K}$ let $\mw_{\mv} \subset \mw$ be the $\mv$-isotypic component of
$\mw$. Then by compactness of $K$ we have $\mw^{fin} = \bigoplus_{\mv \in \hat{K}}
\mw_{\mv}$. Furthermore, if $K$ acts unitarily on $\mw$, then the isotypic components
are orthogonal to each other.

In order to apply this to our setting we consider $\Gamma(\ch_{p,q})$ as a
pre-Hilbert space with respect to the $L^2$-norm, on which the compact group $K$ acts
unitarily. Then its completion $L^2(\ch_{p,q})$ is naturally a $K$-representation for
which the space of smooth sections $\Gamma(\ch_{p,q})$ is dense. Moreover, we know
from \cite{knapp}*{Proposition 8.5} that $K$-finite vectors are smooth. Now for $\mv
\in \hat{K}$ we consider the isotypic component $L^2(\ch_{p,q})_{\mv}$, which
contains $\Gamma(\ch_{p,q})_{\mv} := L^2(\ch_{p,q})_{\mv} \cap \Gamma(\ch_{p,q})$ as
a dense subset. We will refer to $\Gamma(\ch_{p,q})_{\mv}$ as the $\mv$-isotypic
component of $\Gamma(\ch_{p,q})$ by abuse of terminology. In order to analyse the
image of $\uPhi_{p,q}$ we first restrict to an arbitary isotypic component of
$\Gamma(\ch_{p,q})$, which in our setting turns out to be irreducible. More
explicitly, we have

\begin{prop}\label{prop_mult_1}
 Let $G = \SU(n+1,1)$ and $K$ its maximal compact subgroup. Let $\mv$ be a complex
 finite dimensional irreducible $K$-representation. Then the dimension of
 $\Hom_K(\mv, \Gamma(\ch_{p,q}))$ is at most one.
\end{prop}
\begin{proof}
 Follows from \cite{baldoni}*{Theorem 4.4} by Frobenius reciprocity.
\end{proof}

We refer to this result as the multiplicity $1$ property and call $\sigma \in
\Hom_K(\mv, \Gamma(\ch_{p,q}))$ a \emph{Frobenius form}. Note that Proposition
\ref{prop_mult_1} implies that each element in the $\mv$-isotypic component of
$\Gamma(\ch_{p,q})$ can be written as $\sigma(v)$, where the vector $v \in \mv$ is
unique up to multiples. In order to find an explicit formula for the image of
$\uPhi_{p,q}$ on $\Gamma(\ch_{p,q})_{\mv}$ in terms of the Poincar\'{e}-ball model of
hyperbolic space we first analyse the relations between Frobenius forms and the
BGG-operators as well as their adjoints.

\subsection{Frobenius forms and differential operators}
For $0 \le p+q\le n$ and $\mv \in \hat{K}$ let $\sigma \in \Hom_K(\mv,
\Gamma(\ch_{p,q}))$ with $\sigma \neq 0$. Since the partial BGG-operators as well as
their adjoints are both $K$-equivariant differential operators their composition with
$\sigma$ are again Frobenius forms. Moreover, the kernel of $\cd^*$ is closed, so
since $\mv$ is irreducible we either have $\sigma(\mv) \subseteq \ker(\cd^*)$ or
$\sigma(\mv) \cap \ker(\cd^*) = \{0\}$, and similarly for $\ocd^*$. This shows that
for each irreducible $K$-representation the subspace $\sigma(\mv) =
\Gamma(\ch_{p,q})_{\mv}$ of $\Gamma(\ch_{p,q})$ falls into one of four categories,
which will influence the properties of its image under $\uPhi_{p,q}$.

Note that by adjointness the condition $\sigma(\mv) \cap \ker(\cd^*) = \{0\}$ implies
that $\sigma = \cd \hat{\sigma}$ for an element $\hat{\sigma} \in \Hom_K(\mv,
\Gamma(\ch_{p-1,q}))$, which has to be a multiple of $\cd^*\sigma \neq 0$ due to
Proposition \ref{prop_mult_1}. If in addition $\sigma(\mv) \subseteq \ker(\ocd^*)$,
then $\hat{\sigma}(\mv) \subseteq \ker(\cd^*) \cap \ker(\ocd^*)$ since the adjoint
operators anticommute. On the other hand, if additionally $\sigma(\mv) \cap
\ker(\ocd^*) = \{0\}$, then adjointness implies that $\hat{\sigma}$ is a multiple of
$\cd^*\ocd\ocd^*\sigma \neq 0$, which in turn is a nontrivial multiple of $\ocd
\cd^*\ocd^*\sigma$ due to Corollary \ref{cor_commutator_relation_BGG_adjoint}. In
particular, we can write $\sigma = \cd\ocd\tilde{\sigma}$ with $\tilde{\sigma} \in
\Hom_K(\mv, \Gamma(\ch_{p-1,q-1}))$ satisfying $\tilde{\sigma}(\mv) \subseteq
\ker(\cd^*) \cap \ker(\ocd^*)$. By exchanging the r\^{o}les of $\cd^*$ and $\ocd^*$
we obtain analogous results, implying that the case of $\sigma(\mv)$ being contained
in the intersection of the kernels of $\cd^*$ and $\ocd^*$ induces all other cases by
taking images of $\sigma$ under partial BGG-operators.

Next, Proposition \ref{prop_mult_1} implies that each $K$-invariant differential operator $A \colon \Gamma(\ch_{p,q}) \to \Gamma(\ch_{p,q})$ satisfies $A \circ \sigma = \mu(A)\sigma$ with eigenvalue $\mu(A) \in \mc$. Now using the identification of $\ch_{p,q}$ with $\Lambda^{p,q}_0 H^*$ we know from Proposition \ref{prop_operators_pseudo_hermitian}(i) that $\cl_{\zeta}$ induces an invariant operator on $\Gamma(\ch_{p,q})$, and since it is skew-Hermitian with respect to the $L^2$-product it follows that $\mu(\cl_{\zeta})$ is purely imaginary. Similarly, part (v) of the same Proposition shows that also $\Box_H$, $\oBox_H$ and $\Delta_H$ are invariant differential operators on the homology bundles. Furthermore, by selfadjointness all their eigenvalues are positive real numbers, and by Proposition \ref{prop_operators_pseudo_hermitian}(v) we also have
\begin{align}\label{eq:inequality_eigenvalues}
  \mu(\Delta_H) \pm i(n-(p+q))\mu(\cl_{\zeta}) \ge 0
  \end{align}
  for all $0 \le p+q \le n$. 
\begin{prop}\label{prop_Frobenius_forms_diff_operators}
 Let $0 \le p, q \le n$ with $0 \le p+q \le n$. For an irreducible $K$-representation $\mv$ let $\sigma\in \Hom_K(\mv, \Gamma(\ch_{p,q}))$. Then under the identification of $\ch_{p,q}$ with $\Lambda^{p,q}_0 H^*$ we have:
 \begin{enumerate}[(i)]
  \item Let $\sigma(\mv) \cap \ker(\cd^*) = \{0\}$ and $\sigma(\mv) \subseteq \ker(\ocd^*)$ and write $\sigma = \cd \hat{\sigma}$. Then the $\cl_{\zeta}$-eigenvalues $\lambda$ and $\hat{\lambda}$ of $\sigma$ and $\hat{\sigma}$ coincide, whereas their respective $\Delta_H$-eigenvalues $\mu$ and $\hat{\mu}$ are related via $\mu = \hat{\mu} - i\lambda$. 
  
  \item Let $\sigma(\mv) \subseteq \ker(\cd^*)$ and  $\sigma(\mv) \cap \ker(\ocd^*) = \{0\}$ and write $\sigma = \ocd \hat{\sigma}$. Then the $\cl_{\zeta}$-eigenvalues $\lambda$ and $\hat{\lambda}$ of $\sigma$ and $\hat{\sigma}$ coincide, whereas their respective $\Delta_H$-eigenvalues $\mu$ and $\hat{\mu}$ are related via $\mu = \hat{\mu} + i\lambda$.  

  \item Let $p+q = n$. Then $\sigma(\mv) \subseteq \ker(\cd^*)$ if and only if $\opartial_H \circ \sigma = 0$, and similarly for their complex conjugates. In particular, if $\sigma(\mv) \subseteq \ker(\cd^*) \cap \ker(\ocd^*)$, then the eigenvalues of $\Box_H$, $\oBox_H$ and $\Delta_H$ on $\sigma$ are zero.
 \end{enumerate}
\end{prop}
\begin{proof}
 (i) By the multiplicity $1$ property it follows that $\hat{\sigma}$ is a multiple of $\cd^*\sigma$, which under the identification of $\ch_{p,q}$ with $\Lambda^{p,q}_0 H^*$ corresponds to $\partial_H^*\sigma$. Now both $\cl_{\zeta}$ and $\Box_H$ commute with $\partial_H^*$, and expressing $\Delta_H$ in terms these operators via Proposition \ref{prop_operators_pseudo_hermitian}(v) yields the claim.
 
 (ii) Analogously to (i).

 (iii) If $p+q = n$ then the relations in Proposition
   \ref{prop_operators_pseudo_hermitian}(iv) show that $\partial_H^*\sigma = 0$ is
   equivalent to $0 = [L_H^*, \opartial_H]\sigma = L_H^*\opartial_H\sigma$, and since
   $L_H^*$ is injective on $(n+1)$-forms this is equivalent to
   $\opartial_H\sigma = 0$. In particular, if $\sigma(\mv)$ is contained in the
   kernels of both $\cd^*$ and $\ocd^*$, then $\sigma$ lies in the kernels of
   $\partial_H$, $\opartial_H$, and their adjoints, and hence also in the kernels of
   all Laplace operators.
\end{proof}

Consider the $(n+1)$-st BGG-operator $D_{n+1} \colon \Gamma(\ch_n) \to \Gamma(\ch_{n+1})$. In Theorem \ref{thm_BGG_operator}(ii) we have seen that the image of $\iota_{\zeta} D_{n+1}$ consists of primitive forms, however, the operator does not preserve the $(p,q)$-type. For $p+q = n$ denote by $\pi_{p,q} \colon \ch_n \to \ch_{p,q}$ the canonical projection, then $\underline{D}_{n+1} := \pi_{p,q} \circ \iota_{\zeta} \circ D_{n+1} \colon \Gamma(\ch_{p,q}) \to \Gamma(\ch_{p,q})$ is a $K$-invariant differential operator. 

\begin{cor}\label{cor_D_n_formula}
Let $\mv$ be an irreducible $K$-module and $\sigma \in \Hom_K(\mv, \Gamma(\ch_{p,q}))$ with $p+q = n$. Let $\lambda$, $\mu$ and $\nu$ be the eigenvalues of $\cl_{\zeta}$, $\Delta_H$ and $\underline{D}_{n+1}$ on $\sigma$.
\begin{enumerate}[(i)]
  \item If $\sigma(\mv) \subseteq \ker(\cd^*) \cap \ker(\ocd^*)$, then $\nu = \lambda$.
  \item If $\sigma(\mv) \subseteq \ker(\cd^*)$ and $\sigma(\mv) \cap \ker(\ocd^*) = \{0\}$ then $\nu = \lambda - i\frac{\mu}{2}$.
  \item If $\sigma(\mv) \cap \ker(\cd^*) = \{0\}$ and $\sigma(\mv) \subseteq \ker(\ocd^*)$ then $\nu = \lambda + i\frac{\mu}{2}$.
  \item If $\sigma(\mv) \cap \ker(\cd^*) = \{0\}$ and $\sigma(\mv) \cap \ker(\ocd^*) = \{0\}$, then $\sigma(\mv) \subseteq \ker(D_{n+1})$.
 \end{enumerate}
\end{cor}
\begin{proof} Combining the formula for $D_{n+1}$ in Theorem \ref{thm_BGG_operator}(ii) as well as Proposition \ref{prop_Frobenius_forms_diff_operators}(iii) a direct computation yields (i) - (iii). For the last part note that $D_{n+1} \circ \cd = - D_{n+1} \circ \ocd$ since the BGG-operators form a complex. Thus, writing $\sigma = \cd\ocd\hat{\sigma}$ with $\hat{\sigma} \in \Hom_K(\mv,\Gamma(\ch_{p-1,q-1}))$ and using that $\ocd$ squares to zero we immediately obtain $D_{n+1}\sigma = 0$.
\end{proof}

\subsection{The image of $\uPhi_{p,q}$ on isotypic components}\label{sec_im_iso}
Using the relations between differential operators and Frobenius forms from the
previous section we can now determine an explicit formula for the image of
$\uPhi_{p,q}$ on isotypic components of $\Gamma(\ch_{p,q})$ in terms of the
Poincar\'{e} ball model, where we restrict to the case of $0 \le p, q \le n$ with
$p+q \le n$. In order to do so, recall from Theorem \ref{thm_BGG_operator} that
$\uPhi_{p,q} \circ \cd$ and $\uPhi_{p,q} \circ \ocd$ are nontrivial multiples of
$\partial \circ \uPhi_{p-1,q}$ and $\opartial \circ \uPhi_{p,q-1}$,
respectively. Thus, the discussion at the beginning of the previous section implies
that a formula for $\underline{\Phi}_{p,q}(\alpha)$ in the case
$\alpha \in \ker(\cd^*) \cap \ker(\ocd^*)$ yields formulae for other cases by
applying the operators $\partial$ and $\opartial$.

To continue further define for all integers $k$ with $0 \le k \le n$, all
$\lambda \in i\mr$ and $\mu \in \mr_{\ge 0}$ the function
$F_{k,\lambda,\mu} \colon (0,1) \to \mr$ via
\begin{align*}
  F_{k,\lambda,\mu}(r) :=  r^{2\kappa} {}_2 F_1\left(\kappa + \frac{\lambda}{2i}, \kappa - \frac{\lambda}{2i}, 2\kappa + n+1-k;r^2\right),
 \end{align*}
 with the parameter
 $\kappa := \frac{1}{2}\left(\sqrt{(n-k)^2 + (2\mu + |\lambda|^2)} -(n-k)\right) \in
 \mr_{\ge 0}$. The analytic properties of the hypergeometric function in Proposition
 \ref{prop_hypergeometric}(iii) immediately imply that the limit for $r \to 1$ of the
 $\ell$-th derivative of $F_{k,\lambda, \mu}(r)$ exists for all $k \le n-\ell$ and is
 nonzero. Furthermore, part (iv) of the same Proposition shows that the $\ell$-th
 derivative of $F_{n+1-\ell, \lambda, \mu}$ is contained in $L^1([0,1])$. For later
 use, the explicit formula for the limit of a hypergeometric function from
 Proposition \ref{prop_hypergeometric}(iii) yields that the values of
 $F_{k,\lambda,\mu}$ and its derivative at $r = 1$ are related via
\begin{align}\label{eq:limit_derivative_F}
  F_{k,\lambda,\mu}'(1) = F_{k, \lambda, \mu}(1)(\mu - 2\kappa(n-1-k)).
 \end{align}
for all $\lambda \in i\mr$ and $\mu \in \mr_+$.

\begin{prop}\label{prop_Poisson_function_general}
 Let $0 \le p, q \le n$ with $k := p+q \le n$. Let $\mv \in \hat{K}$ and $\sigma \in \Hom_K(\mv, \Gamma(\ch_{p,q}))$ with $\sigma(\mv) \subseteq \ker(\cd^*) \cap \ker(\ocd^*)$. Let $\alpha \in \Gamma(\ch_{p,q})_{\mv}$ and let $\lambda \in i\mr$ and $\mu \in \mr_{\ge 0}$ the eigenvalues of $\cl_{\zeta}$ and $\Delta_H$ on $\alpha$. Then
 \begin{align*}
  \uPhi_{p,q}(\alpha)|_{B\setminus \{0\}} = \gamma_{p,q} F_{k,\lambda,\mu}(1)^{-1} F_{k,\lambda,\mu}(r) \chi^*\alpha
 \end{align*}
 for a constant $\gamma_{p,q} \in \mc$ independent of $\mv$. In particular,
 $\uPhi_{p,q}(\alpha)$ always extends continuously to the boundary. 
\end{prop}
\begin{proof}
  Using the $1$-form $\tau$ and its complex conjugate as well as the polar
  decomposition on $B$ we can write
 \begin{align*}
 (\uPhi_{p,q} \circ \sigma)|_{B\setminus\{0\}} = \chi^*\sigma_{p,q}(r) + \tau \wedge \chi^*\sigma_{p-1,q}(r) + \otau \wedge \chi^*\sigma_{p,q-1}(r) + \tau \wedge \otau \wedge \chi^*\sigma_{p-1,q-1}(r),
\end{align*}
where $\sigma_{s,t}(r)$ is an element in $\Hom_K(\mv, \Gamma(\Lambda^{s,t}H^*))$, and
our aim is to show that only the first summand contributes nontrivially. In order to
do so, recall from Theorem \ref{thm_BGG_operator}(i) that the image of $\uPhi_{p,q}$
consists of primitive differential forms on $B$. Applying the formulae for $L^*$ from
Proposition \ref{prop_polar_differential} we see that $L^*\circ\uPhi_{p,q} = 0$
implies that the equivariant maps $\sigma_{p,q-1}(r)$, $\sigma_{p-1,q}(r)$ and
$\sigma_{p-1,q-1}(r)$ are in the kernel of $L_H^*$. Thus, by Proposition
\ref{prop_mult_1} they have to be a multiple of $\cd^*\sigma$, $\ocd^*\sigma$ and
$\cd^*\ocd^*\sigma$, respectively, and therefore are trivial by assumption. Thus,
only the first summand is nontrivial, and since $L_H^*\sigma_{p,q}(r) = 0$ we obtain
that, on $B\setminus\{0\}$, $\uPhi_{p,q}\circ\sigma = f(r) \chi^*\sigma$ for a smooth
function $f \in C^{\infty}((0,1))$ again by Proposition \ref{prop_mult_1}.
 
Next, we know from Theorem \ref{thm_BGG_operator} that the image of $\uPhi_{p,q}$
consists of coclosed and harmonic differential forms, and we use these properties to
obtain a differential equation for $f$. Note that by Corollary
\ref{lem_kaehler_equivalence_harmonic_partial} this is equivalent to
$L^*d\circ \uPhi_{p,q} = 0$ and $L^*\opartial\partial\circ\uPhi_{p,q} = 0$. Combining
the formulae for the operators on the pseudo-Hermitian manifold $K/M$ in Proposition
\ref{prop_operators_pseudo_hermitian} with the formulae for the operators on $B$ from
Proposition \ref{prop_polar_differential} a few lines of computation show that
$L^*d\circ\uPhi_{p,q}\circ\sigma$ is trivial independently of the function $f$,
whereas the condition $L^*\opartial\partial\circ\uPhi_{p,q}\circ\sigma = 0$ is
equivalent to the differential equation
\begin{align*}
 f'' + \left(\frac{2(n-k)}{r(1-r^2)} + \frac{1}{r}\right)f' + \left(\frac{\lambda^2}{r^2} - \frac{2\mu}{r^2(1-r^2)}\right)f = 0.
\end{align*}
This is a hypergeometric type equation with two linear independent solutions, and
$F_{k,\lambda,\mu}$ is the bounded solution for $r \to 0$ which is unique up to
multiples. Therefore, we know that $f$ coincides with $F_{k,\lambda,\mu}$ up to a
constant $\gamma_{\mv} \in \mc$ which initially depends on $\mv$. In particular, we have
$\lim_{r \to 1} f(r) =\gamma_{\mv} F_{k,\lambda,\mu}(1)$.

Finally, Frobenius reciprocity (c.f. \cite{cap_slovak}*{Theorem 1.1.4}) shows that
$\Hom_K(\mv, \Gamma(\ch_{p,q}))$ is isomorphic to $\Hom_M(\mv, (\ch_{p,q})_{o})$,
where $(\ch_{p,q})_{o}$ denotes the fibre of the vector bundle $\ch_{p,q}$ over the
origin $o = eM \in K/M$. Since $(\ch_{p,q})_{o}$ is an irreducible $M$-module we
deduce that for $\sigma_{\mv} := \sigma \neq 0$ the induced $M$-equivariant map
$\sigma_{\mv}(o) \colon \mv \to (\ch_{p,q})_{o}$ has to be surjective. Now if $\mw$
is another irreducible $K$-representation with nontrivial
$\sigma_{\mw} \in \Hom_K(\mw, \Gamma(\ch_{p,q}))$ satisfying the same assumptions,
then we can find $v \in \mv$ and $w \in \mw$ so that
$\sigma_{\mv}(v)(o) = \sigma_{\mw}(w)(o)$. But we know from Theorem
\ref{thm_main_zero_boundary_value} that the limits of $\uPhi(\sigma_{\mv}(v))$ and
$\uPhi(\sigma_{\mw}(w))$ coincide along any smooth curve on $\overline{B}$ emanating
from $o$, implying that $\gamma_{p,q} := \gamma_{\mv}F_{k, \lambda, \mu}(1)$ is independent of $\mv$.
\end{proof}
\begin{rem}\label{rem-cont-smooth} In Proposition \ref{prop_L_2_norm}
    below, we will show by a more careful study of analytical properties of the
    involved hypergeometric functions that for $p+q=n$ one actually obtains smooth extensions to
    the boundary. One could try to extend these methods to lower degrees, but since
    such results are not needed for the purpose of the current paper, we don't go
    into this. 
\end{rem}

As mentioned above, the previous Proposition yields explicit formulae for the image of each isotypic component $\Gamma(\ch_{p,q})_{\mv}$ under the operator $\uPhi_{p,q}$ in the Poincar\'{e} ball model by applying the commutation relation $c_{p,q} \partial \circ \uPhi_{p-1,q} = (n+1-(p+q)) \uPhi_{p,q} \circ \cd$ as well as its complex conjugate, c.f. Proposition \ref{prop_Poisson_complex_case_p_q}(iii). Moreover, since we are mostly interested in the properties of the image of
$d\uPhi_n$ we do not consider the
isotypic components whose intersection with both $\ker(\cd^*)$ and $\ker(\ocd^*)$ is
trivial, as these are contained in the kernel of $d\uPhi_{p,q}$. In the other cases, the obtained formulae will include one of the constants $\gamma_{p,q}$, $\gamma_{p-1,q}$ or $\gamma_{p,q-1}$ from Proposition \ref{prop_Poisson_function_general}, depending on the $K$-representation $\mv \in \hat{K}$. However, in the end we want to consider the image of $\uPhi_{p,q}$ for a general section of $\ch_{p,q}$, so we need to find a relation between these coefficients. In order to do so, by Theorem \ref{thm_main_zero_boundary_value} we can compare the image of different isotypic components under $\uPhi_{p,q}$ via their limit towards $o = eM \in K/M$. In particular, this yields a relation between the constants $\gamma_{p,q}$ and $c_{p,q}$ for different values of $p$ and $q$, which in turn shows that $\gamma_{p,q}$ coincides with the inverse of $\lambda_{p,q}$ from Theorem \ref{thm_Poisson_complex_space} up to an overall constant $c \in \mc$ independent of $(p,q)$.

\begin{thm}\label{thm_image_poisson_cases}
 Let $0 \le p, q \le n$ with $p+q =: k \le n$ and let $\mv \in \hat{K}$ with $\sigma \in
 \Hom_K(\mv, \Gamma(\ch_{p,q}))$. Let $\alpha \in \Gamma(\ch_{p,q})_{\mv}$ and denote
 by $\lambda$ and $\mu$ the eigenvalues of $\cl_{\zeta}$ and $\Delta_H$ on $\alpha$,
 respectively. Let $\chi \colon B \setminus \{0\} \to K/M$ be the
 angular map, $\gamma_{p,q} \in \mc$ the constants from Proposition \ref{prop_Poisson_function_general} and $c_{p,q}$ the complex numbers from Proposition \ref{prop_Poisson_complex_case_p_q}.
 \begin{enumerate}[(i)]
 \item If $p+q = n$ and $\sigma(\mv) \subseteq \ker(\cd^*)\cap \ker(\ocd^*)$, then $\uPhi_{p,q}(\alpha)|_{B \setminus \{0\}} = \gamma_{p,q} r^{|\lambda|}
   \chi^*\alpha$. Moreover, we have $\opartial\uPhi_{p,q}(\alpha) = 0$ for
   $\lambda \in i\mr_+$ and $\partial\uPhi_{p,q}(\alpha) = 0$ for
   $\lambda \in -i\mr_+$.
  \item  Let $\sigma(\mv) \cap \ker(\cd^*) = \{0\}$ and $\sigma(\mv) \subseteq \ker(\ocd^*)$ and let $\rho$ be the eigenvalue of $\partial_H\partial_H^*$ on $\alpha$. Then with $f(r) := F_{k-1,\lambda, \mu+i\lambda}(r)$ we have
  \begin{align*}
   \uPhi_{p,q}(\alpha)|_{B \setminus \{0\}}= \frac{c_{p,q}\gamma_{p-1,q}}{n+1-k} f(1)^{-1}\left(f(r) \chi^*\alpha + \rho^{-1}\left(\frac{f'(r)}{2r} + \frac{\lambda f(r)}{2ir^2}\right) \tau \wedge \chi^*\partial_H^*\alpha \right). 
  \end{align*}
   \item For $\sigma(\mv) \in \ker(\cd^*)$ and $\sigma(\mv) \cap \ker(\ocd^*) = \{0\}$ let $\rho$ be the eigenvalue of $\opartial_H\opartial_H^*$ on $\alpha$. Then with $f(r) := F_{k-1,\lambda, \mu-i\lambda}(r)$ we have
  \begin{align*}
   \uPhi_{p,q}(\alpha)|_{B \setminus \{0\}}= \frac{\overline{c_{q,p}}\gamma_{p,q-1}}{n+1-k} f(1)^{-1}\left(f(r) \chi^*\alpha + \rho^{-1}\left(\frac{f'(r)}{2r} - \frac{\lambda f(r)}{2ir^2}\right) \otau \wedge \chi^*\opartial_H^*\alpha \right).
  \end{align*}
 \end{enumerate}
 In particular, $\gamma_{p,q} = c \lambda_{p,q}^{-1}$ for an overall constant $c \in \mc$, where $\lambda_{p,q}$ are the coefficients in the definition of $\uPhi_k$ in Theorem \ref{thm_Poisson_complex_space}. 
\end{thm}
\begin{rem} In principle, it would be possible to normalize the constant
    $c$ that shows up in the theorem, for example, by requiring that one obtains the
    classical Poisson transform in degree $0$. Since this is not helpful for the
    further developments, we decided not to normalize, but prove the following
    results for a general choice of $c$.  
\end{rem}
\begin{proof}
 (i) From Proposition \ref{prop_Frobenius_forms_diff_operators}(iv) we know that
   $\mu = 0$ and thus the parameter $\kappa$ in the formula for
   $\uPhi_{p,q}(\alpha)|_{B \setminus \{0\}}$ from Proposition
   \ref{prop_Poisson_function_general} simplifies to $\kappa =
   \frac{|\lambda|}{2}$. In particular, one of the numbers
   $\kappa \pm \frac{\lambda}{2i}$ is zero, and thus the hypergeometric function in
   the same Proposition degenerates to the constant function $1$. Since $p+q = n$ we
   have $\partial_H\sigma = 0$ and $\opartial_H \sigma = 0$ due to Proposition
   \ref{prop_Frobenius_forms_diff_operators}(iv), which implies that
   $d\uPhi_{p,q}(\alpha)$ is a linear combination of $\tau \wedge \chi^*\alpha$ and
   $\otau \wedge \chi^*\alpha$, and computing these coefficients explicitly the
   claims on $\partial\uPhi_{p,q}(\alpha)$ and $\opartial\uPhi_{p,q}(\alpha)$ follow.
  
  (ii) By adjointness we can write $\sigma = \cd \hat{\sigma}$ with the
    $K$-homomorphism $\hat{\sigma} \in \Hom_K(\mv, \Gamma(\ch_{p-1,q}))$ satisfying
    $\hat{\sigma}(\mv) \subseteq \ker(\cd^*) \cap \ker(\ocd^*)$. Since
    $\alpha = \sigma(v)$ for $v \in \mv$ we thus obtain $\alpha = \cd\hat{\alpha}$
    with $\hat{\alpha} = \hat{\sigma}(v) \in \Gamma(\ch_{p-1,q})_{\mv}$. Applying
    Proposition \ref{prop_Poisson_complex_case_p_q}(iii) it follows that
    $\uPhi_{p,q}(\alpha)$ coincides with 
    $\frac{c_{p,q}}{n+1-k} \partial\uPhi_{p-1,q}(\hat{\alpha})$ and therefore
  \begin{align} \label{eq:proof_formula_poisson_case_2}
   \uPhi_{p,q}(\alpha)|_{B \setminus \{0\}}= \frac{c_{p,q}\gamma_{p,q}}{n+1-k} F_{n-1,\hat{\lambda}, \hat{\mu}}(1)^{-1} \partial(F_{n-1,\hat{\lambda}, \hat{\mu}}(r)\chi^*\hat{\alpha})
  \end{align}
due to Proposition \ref{prop_Poisson_function_general}, where $\hat{\lambda}$ and $\hat{\mu}$ are the eigenvalues of $\cl_{\zeta}$ and $\Delta_H$ on $\hat{\alpha}$, respectively. From Proposition \ref{prop_Frobenius_forms_diff_operators}(iii) we know that $\hat{\lambda} = \lambda$ and $\hat{\mu} = \mu + i\lambda$, and the differential form $\hat{\alpha}$ can be explicitly written as $\rho^{-1}\partial_H^*\alpha$. Expansion of the right hand side of \eqref{eq:proof_formula_poisson_case_2} using part (ii) of Proposition \ref{prop_polar_differential} yields the claimed formula.

(iii) Since the map $\uPhi_{p,q}$ is compatible with complex conjugation due to Proposition \ref{prop_Poisson_complex_case_p_q}(ii) the claim follows from (ii) due to $\cl_{\zeta}$ acting on $\overline{\alpha}$ by $\overline{\lambda} = -\lambda$.

For the last part assume that $\sigma_{\mv} \in \Hom_K(\mv, \Gamma(\ch_{p,q}))$ satisfies $\sigma_{\mv}(\mv) \subseteq \ker(\cd^*) \cap \ker(\ocd^*)$. Moreover, let $\mw \in \hat{K}$ with $\sigma_{\mw} \in \Hom_K(\mw, \Gamma(\ch_{p,q}))$ as in (ii). As in the proof of Proposition \ref{prop_Poisson_function_general} Frobenius reciprocity implies that we can find $v\in \mv$ and $w \in \mw$ so that $\sigma_{\mv}(v)(o) = \sigma_{\mw}(w)(o)$ with $o := eM \in K/M$. Furthermore, we know from Theorem \ref{thm_main_zero_boundary_value} that the limits of $\uPhi_{p,q}(\sigma_{\mv}(v))$ and $\uPhi_{p,q}(\sigma_{\mw}(w))$ as differential forms on $\partial \overline{B} \cong K/M$  coindice along any smooth curve on $\overline{B}$ emanating from $o$. Under the identification of $\ch_{p,q}$ with $\Lambda^{p,q}_0 H^*$ these limits are given by $\gamma_{p,q} \sigma_{\mv}(v)(o)$ using the formula from Proposition \ref{prop_Poisson_function_general} and $\frac{c_{p,q}\gamma_{p-1,q}}{n+1-k} \sigma_{\mw}(w)(o)$ via the formula in (ii), respectively. Arguing similarly for $\mw \in \hat{K}$ satisfying the assumptions in (iii) it follows that the coeffients $\gamma_{p,q}$ satisfy the relations
\begin{align*}
 (n+1-(p+q))\gamma_{p,q} = \gamma_{p-1,q}c_{p,q} = \gamma_{p,q-1} \overline{c_{q-1,p}},
\end{align*}
which in turn characterize these numbers up to an overall constant. But \eqref{eq:relation_lambda} shows that $\lambda_{p,q}^{-1}$ also satisfy these relations, from which the claim follows.
\end{proof}
 
The previous Theorem also yields a description of the decomposition of the image of $\uPhi_{p,q}$ for $p+q = n$ into $(p,q)$-types, which we state explicitly for later use.

\begin{cor}\label{cor_Poisson_derivative}
For $0 \le p, q \le n$ with $p+q = n$ let $\cl_{\pm} \subseteq \Gamma(\ch_{p,q})$ be the eigenspaces of $\cl_{\zeta}$ with eigenvalues in $\pm i\mr$ and define the subspaces $\cK$, $\cI_1$ and $\cI_2$ of $\Gamma(\ch_{p,q})$ via
\begin{align*}
 \cK &:= \ker(\cd^*)^{\perp} \cap \ker(\ocd^*)^{\perp}, \\
 \cI_1 &:= \left(\ker(\cd^*) \cap \ker(\ocd^*)^{\perp}\right) \oplus \left(\ker(\cd^*) \cap \ker(\ocd^*) \cap \cl_+ \right), \\
 \cI_2 &:= \left(\ker(\cd^*)^{\perp} \cap \ker(\ocd^*)\right) \oplus \left(\ker(\cd^*) \cap \ker(\ocd^*) \cap \cl_- \right).
\end{align*}
For $\mv \in \hat{K}$ and a $K$-invariant subspace $A \subset \Gamma(\ch_{p,q})$ define $A^{\mv} := A \cap \Gamma(\ch_{p,q})_{\mv}$. Then 
\begin{align*}
 d\uPhi_{p,q}|_{\cK^{\mv}} &= 0, & d\uPhi_{p,q}|_{\cI_1^{\mv}} &\subseteq \Omega^{p+1,q}(G/K), & d\uPhi_{p,q}|_{\cI_2^{\mv}} \subseteq \Omega^{p,q+1}(G/K)
\end{align*}
for all irreducible $K$-representations $\mv$.
\end{cor}

\begin{rem}
  By adjointness it follows that $\Gamma(\ch_{p,q})$ decomposes naturally into the
  direct sum of the kernel of $\cd^*$ and its orthogonal complement, which are both
  closed subspaces, and similarly for its complex conjugated operator. Moreover, each
  $K$-finite subrepresentation of $\Gamma(\ch_{p,q})$ decomposes into irreducibles,
  each of which must be contained in one of the subspaces $\cK$, $\cI_1$ and $\cI_2$
  due to irreducibility. This shows that Theorem \ref{thm_image_poisson_cases} yields
  a full description of $d\uPhi_{p,q}(\Gamma(\ch_{p,q}))$ since $K$-finite vectors
  form a dense subspace.
\end{rem}

\subsection{The $L^2$-norm of $d\uPhi_{p,q}$ on isotypic components}\label{sec_L_2_iso}
In the next step we use the explicit formulae for the image of $\uPhi_{p,q}$ on isotypic components from Theorem \ref{thm_image_poisson_cases} to prove that $d\uPhi_{p,q}(\alpha)$ has finite $L^2$-norm for all $\alpha \in \Gamma(\ch_{p,q})_{\mv}$ with $p+q = n$. We do this by first showing that $\uPhi_{p,q}(\alpha) \in \Omega^{p,q}(G/K)$ as well as its exterior derivative are sufficiently bounded on $\overline{B}$, which allows the application of Stokes' Theorem.

\begin{prop}\label{prop_L_2_norm}
  For all $\alpha \in \Gamma(\ch_{p,q})_{\mv}$ with $p+q=n$, the differential form
  $\uPhi_{p,q}(\alpha) \in \Omega^{n}(B, \mc)$ is smooth up to the boundary, whereas
  $d\uPhi_{p,q}(\alpha) = \phi_1 + f(r)\tau \wedge\otau \wedge \phi_2$ for
  $f \in L^1([0,1])$ and $\phi_1, \phi_2 \in \Omega^*(\overline{B}, \mc)$. Moreover,
  the $L^2$-inner product
  $\langle\!\langle d\uPhi_{p,q}(\alpha), d\uPhi_{p,q}(\beta)\rangle\!\rangle$ is
  finite for all $\alpha$, $\beta \in \Gamma(\ch_{p,q})_{\mv}$. In particular,
  $d\uPhi_{p,q}(\alpha)$ is an $L^2$-harmonic differential form.
\end{prop}
\begin{proof}
The pullback of a differential form on $S^{2n+1}$ along the angular map $\chi$ as
well as the differential forms $\tau$ and $\overline{\tau}$ are smoothly defined on
$\overline{B} \setminus \{0\}$. Thus, Theorem \ref{thm_image_poisson_cases} implies
that the restriction to $B \setminus \{0\}$ of both $\uPhi_{p,q}(\alpha)$ and its
exterior derivative are linear combinations of elements of the form $f(r)\phi$, where
$\phi \in \Omega^*(\overline{B}\setminus\{0\}, \mc)$ and $f \colon (0,1) \to \mc$ is
a smooth function in the radial variable whose limit for $r \to 0$ exists. Moreover,
each such $f$ is a linear combination of $F_{k,\lambda,\mu}$ for $k \in \{n, n-1\}$
as well as their first and second derivatives, and the latter can only appear in the
coefficient of an element of the form $\tau\wedge \overline{\tau} \wedge \phi$. Using
the analytic properties of the hypergeometric function from Proposition
\ref{prop_hypergeometric}, we have seen that the functions $F_{k,\lambda,\mu}$ as well
as their first derivative can be smoothly extended to $[0,1]$, whereas
$F_{n-1,\lambda, \mu}'' \in L^1([0,1])$. From this the claims on
$\uPhi_{p,q}(\alpha)$ and $d\uPhi_{p,q}(\alpha)$ follow.

For the computation of the $L^2$-inner product recall from Proposition
\ref{prop_Poisson_complex_case_p_q} that $\uPhi_{p,q}(\alpha)$ and
$\uPhi_{p,q}(\beta)$ are harmonic and coclosed differential forms, which implies
\begin{align}\label{eq:L_2_Poisson_int}
 \langle\!\langle d\uPhi_{p,q}(\alpha), d\uPhi_{p,q}(\beta)\rangle\!\rangle = \int_{B} d\left(\uPhi_{p,q}(\alpha) \wedge \ast d\overline{\uPhi_{p,q}(\beta)}\right).
\end{align}
Since $\ast$ is tensorial, it follows that the integrand of \eqref{eq:L_2_Poisson_int}
can be written as $d\psi_1 + g(r) d(\tau \wedge \otau\wedge \psi_2)$ with $g \in
L^1([0,1])$ and $\psi_1$, $\psi_2 \in \Omega^*(\overline{B}, \mc)$. Now using a
smooth approximation of $g$ with support in $B$ together with Stokes' theorem, an easy
estimate yields that the integral of the second summand is trivial. Thus, we can
apply Stokes' theorem to \eqref{eq:L_2_Poisson_int} and the domain of the resulting
integral is the compact space $\partial B \cong K/M$, showing that the $L^2$-inner
product is finite.
\end{proof}

\begin{rem}
The above realization of $L^2$-harmonic $(n+1)$-forms also includes the construction given in \cite{julg_kasparov}*{Theorem 3.2}. There, the initial data was a differential form $\omega \in \Omega^{p,q}(\mc^{n+1})$ with $p+q = n+1$ and $p$, $q \neq 0$ which is closed and primitive for the Euclidean metric. Furthermore, it was assumed that the coefficients of $\omega$ are homogeneous polynomials in $z$ and $\bz$ of degree $k$ and $m$, respectively. To translate this to our setting we view $\omega$ as an element of the irreducible $K$-representation $\mv := S^{k,m}(\mc)\otimes \Lambda^{p,q}_0\mc^{n+1}$, where $S^{k,m}(\mc)$ is the space of harmonic polynomials in $z$ and $\overline{z}$ of degree $k$ and $m$, respectively. Then $\omega \mapsto \iota_{\overline{T}}\iota_T\omega|_{K/M}$ defines a $K$-homomorphism $\sigma \in \Hom_K(\mv, \Gamma(\ch_{p-1,q-1}))$. Putting $a := k + p$ and $b := m + q$, the formula for $\Delta_H$ in \cite{rumin}*{Proposition 2} yields that $\cl_{\zeta}\sigma = i(a-b)\sigma$ and $\Delta_H\sigma = (a+b+2ab)\sigma$, and inserting this into Proposition \ref{prop_Poisson_function_general} shows that $\uPhi_{p-1,q-1}(\sigma(\omega)) = cF(1)^{-1} r^{a+b}F(t)\chi^*\sigma(\omega)$ with $F(t) := {}_2F_1(a,b,a+b+2;t)$. The differential forms $\opartial(F(t)\sigma(\omega))$ and $-\partial(F(t)\sigma(\omega))$ coincide with the forms $\varphi_T$ and $\varphi_{\overline{T}}$ in \cite{julg_kasparov}, and thus $\uPhi_{p-1,q-1}(\sigma(\omega))$ coincides with the formula in the paper up to normalization with $F(1)^{-1}$ and the factor $r^{a+b}$, the latter appearing due to the pullback along the angular map.
\end{rem}

Although the previous Proposition guarantees that the image of $d\uPhi_{p,q}$ on each
isotypic component of $\Gamma(\ch_{p,q})$ consists of $L^2$-harmonic forms on $G/K$
this is not sufficient to show the same property for the image of $d\uPhi_{p,q}$ on
all of $\Gamma(\ch_{p,q})$. Therefore, we aim to express the $L^2$-norm of
$d\uPhi_{p,q}(\alpha)$ in terms of globally defined operators on $K/M$, which then
can be extended to all of $\Gamma(\ch_{p,q})$. It turns out that a slight variation
of the $L^2$-inner product gives a nicer description in terms of boundary
operators. Explicitly, from \cite{weil}*{Thm 1.3} we know that the Hodge star
operator acts on primitive $(n+1)$-forms on $G/K$ by multiplication with a non-zero 
number, which implies that
\begin{align*}
\int_{G/K} d\uPhi_{p,q}(\alpha) \wedge \overline{d\uPhi_{p,q}(\beta)}
\end{align*}
coincides with the $L^2$-inner product of $d\uPhi_{p,q}(\alpha)$ and $d\uPhi_{p,q}(\beta)$ up to a multiple and is therefore finite. Our aim is to relate this expression to the $(n+1)$-st BGG-operator, for which we use the following formula.

\begin{lemma}\label{lem_L_2_comp}
 Let $0 \le p, q \le n$ with $p+q = n$, let $\mv \in \widehat{K}$ and $\alpha$, $\beta \in \Gamma(\ch_{p,q})_{\mv}$. Then we can write
 \begin{align*}
  \left.\left(\uPhi_{p,q}(\alpha) \wedge d\overline{\uPhi_{p,q}(\beta)}\right)\right|_{B \setminus \{0\}} \equiv (f_1(r) \tau + f_2(r)\otau) \wedge \chi^*(\alpha \wedge \overline{\beta}) \mod \tau \wedge \otau
 \end{align*}
with smooth functions $f_1$, $f_2 \colon (0,1) \to \mc$ whose limits for $r \to 1$ exist. Furthermore, 
\begin{align}\label{eq:L_2_Poisson_reform}
 \int_{G/K} d\uPhi_{p,q}(\alpha) \wedge \overline{d\uPhi_{p,q}(\beta)} = (-1)^ni(f_1(1) - f_2(1)) \int_{K/M}  \alpha \wedge \vartheta \wedge \overline{\beta},
\end{align}
where $\vartheta \in \Omega^1(K/M)$ denotes the contact form.
\end{lemma}

\begin{proof}
 The integrand on the left hand side of $\eqref{eq:L_2_Poisson_reform}$ can be written as the exterior derivative of $\Psi := \uPhi_{p,q}(\alpha) \wedge \overline{d\uPhi_{p,q}(\beta)} \in \Omega^{2n+1}(\overline{B}, \mc)$. Therefore, its wedge product with $\tau \wedge \otau$ has to be trivial, and thus Proposition \ref{prop_L_2_norm} implies that modulo $\tau \wedge \otau$ we can write $\Psi$ as a linear combination of $\tau \wedge \psi_1$ and $\otau \wedge \psi_2$, where $\psi_1$, $\psi_2 \in \Omega^{n,n}(G/K)$ are trivial upon insertion of either vector field $T$, $\overline{T} \in \mathfrak{X}(G/K)$. In particular, for fixed $r \in (0,1)$ we can use the angular map to identify the restriction of $\psi_j$ to the $r$-sphere $S_r$ of $B$ with a section of $\Lambda^{n,n} H^*$. But the latter bundle has rank $1$, implying that we can write $\psi_j = f_j(r) \chi^*(\alpha \wedge \overline{\beta})$. Therefore, $\Psi$ has the claimed form, and the limits of the coefficient functions exist due to Proposition \ref{prop_L_2_norm}. For the integral we write the left hand side of \eqref{eq:L_2_Poisson_reform} as the integral over $d\Psi$, and as in the proof of Proposition \ref{prop_L_2_norm} the component of $\Psi$ containing $\tau \wedge \otau$ does not contribute to this integral. Thus, the formula follows from Stokes' theorem and using that the inclusion $\iota \colon K/M \hookrightarrow \overline{B}$ satisfies $\iota^*\tau = i \vartheta$.
\end{proof}

We will use the Lemma to relate the left hand side of \eqref{eq:L_2_Poisson_reform}
to an integral over $K/M$ involving the $(n+1)$-st BGG-operator. In order to do so,
recall that for $0 \le p, q\le n$ with $p+q = n$ the invariant operator
$\underline{D}_{n+1}$ on $\Gamma(\ch_{p,q})$ is the component of
$\iota_{\zeta}D_{n+1}$ preserving the bidegree. By the multiplicity $1$ property from
Proposition \ref{prop_mult_1} this operator acts by a scalar $\nu \in i\mr$ on each
$\mv$-isotypic component of $\Gamma(\ch_{p,q})$, and we computed this factor
explicitly in Corollary \ref{cor_D_n_formula}.

\begin{thm}\label{thm_relation_integral_D_n}
 For $0 \le p, q \le n$ with $p+q = n$ and $\mv \in \widehat{K}$ let $\alpha$, $\beta
 \in \Gamma(\ch_{p,q})_{\mv}$. Let $\lambda_{p,q} \in \mc$ be the coefficients in the
 definition of $\uPhi_{n}$ from Theorem \ref{thm_Poisson_complex_space} and let
 $c\in\mc$ be the constant from Theorem \ref{thm_image_poisson_cases} (which is
 independent of $\mv$, $p$, and $q$). Then we have
 \begin{align}\label{eq:L_2_Frobenius_global_formula}
  \int_{G/K} d\uPhi_{p,q}(\alpha) \wedge \overline{d\uPhi_{p,q}(\beta)} = |\lambda_{p,q}|^{-2}  |c|^2\int_{G/P} \alpha \wedge D_{n+1}\overline{\beta}.
 \end{align}
Equivalently, if $\nu$ is the eigenvalue of $\underline{D}_{n+1}$ on $\Gamma(\ch_{p,q})_{\mv}$, then 
\begin{align*}
 \langle\!\langle d\uPhi_{p,q}(\alpha), d\uPhi_{p,q}(\beta)\rangle\!\rangle_{L^2} = |\lambda_{p,q}|^{-2} |c|^2 |\nu| \langle\!\langle \alpha, \beta \rangle\!\rangle_{L^2}.
\end{align*}
\end{thm}
\begin{proof}
By comparing bidegrees we can replace $D_{n+1}\overline{\beta}$ by $\vartheta \wedge
\underline{D}_{n+1}\overline{\beta}$ on the right hand side of
\eqref{eq:L_2_Frobenius_global_formula}. Since $\nu \in i\mr$ is purely imaginary
complex conjugation yields $\underline{D}_{n+1}\overline{\beta} = -\nu
\overline{\beta}$. Therefore, if $f_1(r)$ and $f_2(r)$ are as in Lemma
\ref{lem_L_2_comp} it suffices to prove
\begin{align}\label{eq:eigenvalue_D_n_f_i}
 (-1)^n i(f_2(1)-f_1(1)) = |\lambda_{p,q}|^{-2} |c|^2\nu
\end{align}
for the constant $c \in \mc$ coming from the explicit
formulae for the image of $\uPhi_{p,q}$ on isotypic components in Theorem
\ref{thm_image_poisson_cases}. In order to do so, let $\lambda$ and $\mu$ be the
eigenvalues of $\cl_{\zeta}$ and $\Delta_H$ on $\Gamma(\ch_{p,q})_{\mv}$,
respectively, and write $\alpha = \sigma(v)$ and $\beta =\sigma(w)$ for $v$, $w \in
\mv$ and $\sigma \in \Hom_K(\mv, \Gamma(\ch_{p,q}))$, which is unique up to
multiples.
 
 (i) If $\sigma(\mv) \subseteq \im(\cd^*) \cap \im(\ocd^*)$, then since $\uPhi_{p,q}
  \circ \cd$ is a nontrivial multiple of $\partial \circ \uPhi_{p-1,q}$ as well as
  for their complex conjugates due to Theorem \ref{thm_PT_rank_1} we obtain
  $d\uPhi_{p,q}\beta = 0$, whereas $D_{n+1}\overline{\beta} = 0$ due to Corollary
  \ref{cor_D_n_formula}(iv).

(ii) If $\sigma(\mv) \subseteq \ker(\cd^*) \cap \ker(\ocd^*)$ then
  $\uPhi_{p,q}(\alpha) = \lambda_{p,q}^{-1} cr^{|\lambda|}\chi^*\alpha$ and $\nu = \lambda$. Since
  $\partial_H\beta = 0$ and $\opartial_H\beta = 0$ due to Proposition
  \ref{prop_Frobenius_forms_diff_operators}(iv) the exterior derivative of
  $\uPhi_{p,q}(\beta)$ is a linear combination of $\tau \wedge \chi^*\beta$ and
  $\otau \wedge \chi^*\beta$, and determining the coefficients explicitly yields
  \eqref{eq:eigenvalue_D_n_f_i}.

 (iii) If $\sigma(\mv) \cap \ker(\cd^*) = \{0\}$ and $\sigma(\mv) \subset
   \ker(\ocd^*)$, then we have $\uPhi_{p,q}(\alpha) = 2 \lambda_{p,q}^{-1}c\mu^{-1} \partial (f(r)
   \chi^*\partial_H^*\alpha)$ with $f(r) = F_{n-1, \lambda,
     \mu+i\lambda}(1)^{-1}F_{n-1,\lambda, \mu+i\lambda}(r)$, whereas $\nu = \lambda -
   \frac{i}{2}\mu$. Define the functions $g_{\pm}(r) := \frac{1}{2r} (f' \pm i\lambda
   f)$ and note that since $f'(1) = \mu + i\lambda$ due to
   \eqref{eq:limit_derivative_F} we get $g_-(1) = \frac{\mu}{2}$ and $g_+(1) =
   \frac{\mu}{2} + i\lambda$. Moreover, a few lines of computation using the calculus
   on $K/M$ and $B$ from Propositions \ref{prop_operators_pseudo_hermitian} and
   \ref{prop_polar_differential} yields that, modulo $\tau \wedge \otau$, we get 
$$
\opartial\uPhi_{p,q}(\beta) \equiv c\lambda_{p,q}^{-1}\big((f-2\mu^{-1}g_- r^2) \chi^*\opartial_H\beta
- 2\mu^{-1}g_- \tau \wedge \chi^*\opartial_H\partial_H^*\beta +  g_+ \otau \wedge
\chi^*\beta\big)
$$
Since formula \eqref{eq:L_2_Poisson_reform} in Lemma \ref{lem_L_2_comp} only depends
on the values of the functions $f_1$ and $f_2$ at $r = 1$ we can omit all summands in $\opartial\uPhi_{p,q}(\beta)$ whose boundary value is trivial. A direct computation shows that this is the case for the coefficient function of $\chi^*\opartial_H\beta$. Moreover, $\opartial_H\partial_H^*\beta$ has values in $\ch_{p-1,q+1}$ and thus the wedge product of its complex conjugate with $\uPhi_{p,q}(\alpha)$ is trivial. Therefore, we can replace $\opartial\uPhi_{p,q}(\beta)$ with $g_+ \otau \wedge \chi^*\beta$, and a direct computation shows that $f_2(1) = 0$ and $f_1(1) = (-1)^n\lambda_{p,q}^{-1} |c|^2(\frac{\mu}{2} + i \lambda)$, proving \eqref{eq:eigenvalue_D_n_f_i}.

(iv) The case $\sigma(\mv) \subseteq \ker(\cd^*)$ and $\sigma(\mv) \cap \ker(\ocd^*) = \{0\}$ follows from (iii) since $\uPhi_{p,q}$ is compatible with complex conjugation due to Theorem \ref{thm_PT_rank_1}.

 To prove the last statement, we use that the Hodge star operator acts on a primitive
 $(s,t)$-form on $G/K$ with $s+t = n+1$ by multiplication with
 $(-1)^{\frac{(n+1)(n+2)}{2}} i^{s-t}$ due to \cite{weil}*{Thm 1.3}, and a similar
 proof yields that the Hodge star operator $\ast_H$ acts on primitive sections of
 $\Lambda^{p,q} H^*$ with $p+q = n$ by multiplication with $(-1)^{\frac{n(n+1)}{2}}
 i^{p-q}$. Now from Theorem \ref{thm_image_poisson_cases} we know that
 $\uPhi_{p,q}(\beta)$ is in the kernel of $\partial$ or $\opartial$, and rest follows
 from an easy computation.
\end{proof}

\section{Poisson transforms and discrete series representations}\label{sec_DSR}
As a last step, we analyze the real operator $\uPhi_n \colon \Gamma(\ch_n) \to
\Omega^n(G/K)$ defined in Theorem \ref{thm_Poisson_complex_space} to prove that the
image of $d\uPhi_n$ consists of $L^2$-forms. Now it is well known, see
  \cite{pedon_complex} that the space $L^2\Omega^{n+1}_{\Delta}(G/K)$ of
  $L^2$-harmonic $(n+1)$-forms on $G/K$ is a realization of the direct sum of all
  discrete series representations of $G$ with trivial infinitesimal character.
Recall that a discrete series representation for a linear connected semisimple Lie
group $G$ is an irreducible representation on a Hilbert space $\ch$ so that its
matrix coefficients are contained in $L^2(G)$. By \cite{knapp}*{Theorem 12.20} such
representations can only exist if the ranks of $G$ and its maximal compact subgroup
$K$ coincide. Furthermore, if this is the case then the number of inequivalent
discrete series representations with the same infinitesimal character equals
$|W_G|/|W_K|$, where $W_G$ and $W_K$ are the Weyl groups of $G$ and $K$,
respectively. Now in \cite{borel}*{Theorem A} it was shown that discrete series
representations with trivial infinitesimal character can be realized as
$L^2$-harmonic differential forms on $G/K$, which can only exist in degree
$\frac{1}{2} \dim_{\mr}(G/K)$. In our setting we have $|W_G|/|W_K| = n+2$, and it was
computed in \cite{borel_wallach}*{Theorem 4.11} and \cite{pedon_complex}*{Theorem
  3.4} that the inequivalent discrete series representations with trivial
infinitesimal character correspond to the decomposition of $L^2$-harmonic
$(n+1)$-forms into $(p,q)$-types.

\subsection{The image of $d\uPhi_n$ on $K$-finite vectors}
Recall from Theorem \ref{thm_Poisson_complex_space} that the real operator $\uPhi_n
\colon \Gamma(\ch_n) \to \Omega^n(G/K)$ was defined as a linear combination of the
complex operators $\uPhi_{p,q} \colon \Gamma(\ch_{p,q}) \to
\Omega^{p,q}(G/K)$. Moreover, the space of $K$-finite vectors of $\Gamma(\ch_n)$ is a
dense subspace which decomposes into the sum of isotypic components. Thus, in order
to gain information about the kernel and the image of $d\uPhi_n$ we will first
consider its restriction to isotypic components and apply the results from Sections
\ref{sec_im_iso} and \ref{sec_L_2_iso}, and afterwards use a completion argument to
obtain general results.

Given a $K$-finite section $\alpha \in \Gamma(\ch_n)$ we can first decompose it into
$(p,q)$-types $\alpha_{p,q}$ and afterwards into isotypic components
$\alpha_{p,q}^{\mv}$. Since the projection $\ch_n \to \ch_{p,q}$ is $K$-equivariant
and tensorial it follows that $\Gamma(\ch_n)_{\mv} \cap \Gamma(\ch_{p,q})$ coincides
with the $\mv$-isotypic component $\Gamma(\ch_{p,q})_{\mv}$. In particular, we can
write each element in $\Gamma(\ch_n)_{\mv}$ as a sum of elements in
$\Gamma(\ch_{p,q})_{\mv}$. Since different isotypic components and $(p,q)$-types are
orthogonal it follows that $\alpha = 0$ if and only if $\alpha^{\mv}_{p,q} = 0$ for
all $\mv \in \widehat{K}$ and all $0 \le p, q \le n$ with $p+q = n$. Similarly, if
$\beta \in \Gamma(\ch_n)$ is a $K$-finite section, then using $K$-invariance of the
$L^2$-inner product it follows that $\langle\!\langle \alpha, \beta \rangle\!\rangle$
coincides with the sum of $\langle\!\langle \alpha^{\mv}_{p,q}, \beta^{\mv}_{p,q}
\rangle\!\rangle$ over all $\mv \in \hat{K}$ and $p+q = n$.

Next, the space $\Omega^{n+1}(G/K)$ is naturally a unitary $K$-representation, and
since $d\uPhi_n$ is $K$-equivariant it maps $\mv$-isotypic components to
$\mv$-isotypic components. Thus, by linearity of $d\uPhi_n$ it follows that
$d\uPhi_n(\alpha) = 0$ if and only if $d\uPhi_n(\alpha^{\mv}) = 0$ for all $\mv \in
\hat{K}$. Furthermore, for the $L^2$-inner product we obtain that $\langle\!\langle
d\uPhi_n(\alpha), d\uPhi_n(\beta)\rangle\!\rangle$ is the sum of $\langle\!\langle
d\uPhi_n(\alpha^{\mv}) , d\uPhi_n(\beta^{\mv})\rangle\!\rangle$ over all $\mv \in
\hat{K}$. However, since $d\uPhi_n$ does not preserve the bidegree it follows that
differential forms $d\uPhi_n(\alpha_{p,q})$ can be linearly dependent for different
values of $p$ and $q$. In the following Lemma we show that this only happens in one
specific situation.

\begin{lemma}\label{lem_cases_pairs_degree}
 For $\mv \in \hat{K}$ and $1 \le s,t \le n$ with $s+t = n+1$ let $\alpha_{s-1,t} \in
 \Gamma(\ch_{s-1,t})_{\mv}$ and $\alpha_{s,t-1} \in \Gamma(\ch_{s,t-1})_{\mv}$ and
 assume that both $\partial \uPhi_n(\alpha_{s-1,t})$ and $\opartial
 \uPhi_n(\alpha_{s,t-1})$ are nontrivial. Then $\alpha_{s-1,t} \in \ker(\cd^*)$ and
 $\alpha_{s,t-1} \in \ker(\ocd^*)$. Moreover, we have one of two possible cases:
\begin{enumerate}[(i)]
 \item The elements $\alpha_{s-1,t}$ and $\alpha_{s,t-1}$ are contained in the kernels of both $\cd^*$ and $\ocd^*$. Moreover, $\partial \uPhi_n(\alpha_{s-1,t})$ and $\opartial \uPhi_n(\alpha_{s,t-1})$ are orthogonal with respect to the $L^2$-inner product.
 \item We have $\alpha_{s-1,t} \in \im(\ocd)$ as well as $\alpha_{s,t-1} \in \im(\cd)$, and we can write
 \begin{align*}
  \partial \uPhi_n(\alpha_{s-1,t}) + \opartial \uPhi_n(\alpha_{s,t-1}) = \partial\opartial \uPhi_{n-1}(\hat{\alpha} _1- \hat{\alpha}_2)
 \end{align*}
 for $\hat{\alpha}_j \in \Gamma(\ch_{s-1,t-1})_{\mv}$ with $\alpha_{s-1,t} = \ocd \hat{\alpha}_1$ and $\alpha_{s,t-1} = \cd\hat{\alpha}_2$.
\end{enumerate}
\end{lemma}
\begin{proof}
 Let $\sigma_{p,q} \in \Hom_K(\mv, \Gamma(\ch_{p,q}))$ so that $\alpha_{p,q} \in
 \im(\sigma_{p,q})$. By the multiplicity $1$-property in Proposition
 \ref{prop_mult_1} we know that $\sigma_{p,q}$ is unique up to multiples, and by
 assumption both $\sigma_{s-1,t}$ and $\sigma_{s,t-1}$ are nonzero. From Corollary
 \ref{cor_Poisson_derivative} it follows that $\partial \uPhi_n(\alpha_{s-1,t}) \neq
 0$ only if $\alpha_{s-1,t} \in \ker(\cd^*)$, which by irreducibility of $\mv$ is
 equivalent to $\sigma_{s-1,t}(\mv) \subseteq \ker(\cd^*)$, and similarly for
 $\opartial \uPhi_n(\alpha_{s,t-1})$. Furthermore, by irreducibility of $\mv$ it
 follows that $\sigma_{s-1,t}(\mv)$ is either a subspace of $\ker(\ocd^*)$ or
 has trivial intersection with it.
 
 If $\sigma_{s-1,t}(\mv) \cap \ker(\ocd^*) = \{0\}$ then we can write $\sigma_{s-1,t}
 = \ocd\hat{\sigma}$ for $\hat{\sigma} \in \Hom_K(\mv, \Gamma(\ch_{s-1,t-1}))$. Since
 $\alpha_{s-1,t} = \sigma_{s-1,t}(v)$ for $v \in \mv$ this shows that $\alpha_{s-1,t}
 \in \im(\ocd)$. Furthermore $\cd\hat{\sigma} \neq 0$, because otherwise we have
\begin{align*}
 \partial\uPhi_n(\alpha_{s-1,t}) = \partial\uPhi_n(\ocd\hat{\sigma}(v)) = \partial\opartial \uPhi_{n-1}(\hat{\sigma}(v)) = -\opartial\partial \uPhi_{n-1}(\hat{\sigma}(v)) = -\opartial \uPhi_n(\cd \hat{\sigma}(v)) = 0,
\end{align*}
contradicting our assumption. By the multiplicity $1$ property in Proposition
\ref{prop_mult_1}, the map $\cd\hat{\sigma}$ has to be a nontrivial multiple of
$\sigma_{s,t-1}$, which shows that $\alpha_{s,t-1} \in
\im(\cd)$. Moreover, by adjointness this implies $\cd^*\sigma_{s,t-1} \neq
0$ and therefore $\sigma_{s,t-1}(\mv) \cap \ker(\cd^*) = \{0\}$.
Additionally, from Theorem \ref{thm_Poisson_complex_space} we know that $\uPhi_n
\circ D_n = d \circ \uPhi_{n-1}$, which implies the claim for the sum of the
derivatives. Via an analogous argument we also obtain $\sigma_{s-1,t}(\mv) \cap
\ker(\ocd^*) = \{0\}$ is equivalent to $\sigma_{s,t-1}(\mv) \cap \ker(\cd^*) =
\{0\}$, proving that we are always in one of the two claimed cases.

In the other case $\sigma_{s-1,t}(\mv)$ and $\sigma_{s,t-1}(\mv)$ are contained in
both kernels of $\cd^*$ and $\ocd^*$. Moreover, we have seen in the proof of Theorem
\ref{thm_relation_integral_D_n} that $\partial \uPhi_n(\alpha_{s-1,t})$ and
$\opartial \uPhi_n(\alpha_{s,t-1})$ are a multiple of $\tau \wedge
\chi^*\alpha_{s-1,t}$ and $\otau \wedge \chi^*\alpha_{s,t-1}$, respectively, which
are orthogonal with respect to the $L^2$-inner product.
\end{proof}

Now we have all ingredients at hand to prove the main result of this article, which
explicitly relates the BGG-complex to the discrete series representations of $G$.

\begin{thm}\label{prop_L_2_general}
  Let $G = \SU(n+1,1)$, $K \subset G$ its maximal compact subgroup and $P \subset G$
its minimal parabolic subgroup and consider the $G$-equivariant operator $\uPhi_n \colon \Gamma(\ch_n) \to \Omega^n(G/K)$ from Theorem \ref{thm_Poisson_complex_space}.
\begin{enumerate}[(i)]
 \item We have $\ker(d\uPhi_n) = \im(D_{n})$, where $D_{n}$ is the $n$-th BGG operator.
 \item For all $\alpha$, $\beta \in \Gamma(\ch_n)$ we have
\begin{align}\label{eq:int_equality_D_n_general}
 \int_{G/K} d\uPhi_n(\alpha) \wedge d\uPhi_n(\overline{\beta}) = |c|^2\int_{G/P} \alpha \wedge D_{n+1}\overline{\beta}
\end{align}
for the constant $c \in \mc$ from Theorem \ref{thm_image_poisson_cases}. In particular, $\langle\!\langle d\uPhi_n(\alpha),
d\uPhi_n(\beta) \rangle \!\rangle$ is finite, and the differential form
$d\uPhi_n(\alpha) \in \Omega^{n+1}(G/K)$ is $L^2$-harmonic.
\item The Poisson transform $d\uPhi_n$ maps $\Gamma(\ch_n)/\im(D_{n})$
    isomorphically onto a dense subspace of the space $L^2\Omega^{n+1}_{\Delta}(G/K)$
    of $L^2$-harmonic $(n+1)$-forms on $G/K$, which is the direct sum of all discrete
    series representations of $G$ with trivial infinitesimal character.
\end{enumerate}
\end{thm}

\begin{proof}
(i) By Theorem \ref{thm_Poisson_complex_space} we have $\uPhi_n \circ D_n = d
    \circ \uPhi_{n-1}$ and thus $\im(D_n) \subseteq \ker(d\uPhi_n)$. For the other
    direction note that both $\ker(d\uPhi_n)$ and $\im(D_n) = \ker(D_{n+1})$ are
    closed $K$-invariant subspaces of $\Gamma(\ch_n)$. Moreover, the subspace of
    $K$-finite vectors $\ker(d\uPhi_n)^{fin}$ of $\ker(d\uPhi_n)$ is
    dense. Therefore, if $\ker(d\uPhi_n)$ is contained in $\im(D_n)$ then taking
    the closure of the chain of inclusions $\ker(d\uPhi_n)^{fin} \subseteq \im(D_n)
    \subseteq \ker(d\uPhi_n)$ implies the claim. Furthermore, we have seen that for a
    $K$-finite section $\alpha \in \Gamma(\ch_n)$ we have $d\uPhi_n(\alpha) = 0$ if
    and only if $d\uPhi_n$ is trivial on each isotypic component of $\alpha$.
  
  Therefore, for a fixed $\mv \in \hat{K}$, let $\alpha \in \Gamma(\ch_n)_{\mv}$ with
  $d\uPhi_n(\alpha) = 0$ and let $\alpha = \sum_{p+q=n} \alpha_{p,q}$ be its
  decomposition into $(p,q)$-types. From Corollary \ref{cor_Poisson_derivative} it
  follows that $\uPhi_n(\alpha_{p,q})$ is contained in $\ker(\partial)$ or in
  $\ker(\opartial)$, and $d\uPhi_n(\alpha_{p,q}) = 0$ if and only if $\alpha_{p,q}
  \in \ker(\cd^*)^{\perp} \cap \ker(\ocd^*)^{\perp}$. Writing $\alpha_{p,q} =
  \sigma_{p,q}(v)$ with $\sigma_{p,q} \in \Hom_K(\mv, \Gamma(\ch_{p,q}))$ this means
  that $\sigma_{p,q}(\mv) \cap \ker(\cd^*) = \{0\}$ and $\sigma_{p,q}(\mv) \cap
  \ker(\ocd^*) = \{0\}$, and by Corollary \ref{cor_D_n_formula} this implies that
  $D_{n+1} \circ \sigma_{p,q} = 0$ and thus $\alpha_{p,q} \in \ker(D_{n+1}) =
  \im(D_n)$.
  
  By symmetry we can assume w.l.o.g that $\partial\uPhi_n(\alpha_{s-1,t}) \neq 0$ for
  indices $0 \le s, t \le n+1$ with $s+t = n+1$. Since $d\uPhi_n(\alpha) = 0$, this
  implies $\partial \uPhi_n(\alpha_{s-1,t}) = -\opartial \uPhi_n(\alpha_{s,t-1})$ by
  comparing $(p,q)$-types. From Lemma \ref{lem_cases_pairs_degree} we know that this
  is only possible if $\alpha_{s-1,t} = \ocd \hat{\alpha}_1$ and $\alpha_{s,t-1} =
  \cd \hat{\alpha}_2$ for $\hat{\alpha}_j \in \Gamma(\ch_{s-1,t-1})_{\mv}$, and we
  obtain the relation $\partial\opartial\uPhi_{n-1}(\hat{\alpha}_1) =
  \partial\opartial\uPhi_{n-1}(\hat{\alpha}_2)$. From the explicit formulae for the
  image of the Poisson transform on the Poincar\'{e} ball model in Theorem
  \ref{thm_image_poisson_cases} this is only possible if $\hat{\alpha}_1 =
  \hat{\alpha}_2$ and thus $\alpha_{s-1,t} + \alpha_{s,t-1} \in \im(D_n)$.
  
  (ii) As in the proof of Theorem \ref{thm_relation_integral_D_n} we obtain that the
    left hand side of \eqref{eq:int_equality_D_n_general} is a nontrivial multiple of
    the $L^2$-inner product on $G/K$, whereas the right hand side can be written as a
    multiple of $\langle\!\langle \alpha,
    \iota_{\zeta}D_{n+1}\beta\rangle\!\rangle$. Now assume that we have shown the
    claimed equation for all $K$-finite vectors $\alpha$, $\beta \in
    \Gamma(\ch_n)$. Then for fixed $\beta$ both sides of
    \eqref{eq:int_equality_D_n_general} are continuous functions $\Gamma(\ch_n) \to
    \mr$ which coincide on the dense subspace of $K$-finite vectors and thus have to
    be equal on all of $\Gamma(\ch_n)$. Now for fixed $\alpha \in \Gamma(\ch_n)$ the
    same argument shows equality of \eqref{eq:int_equality_D_n_general} for all
    $\beta \in \Gamma(\ch_n)$. Moreover, isotypic components of $\Gamma(\ch_n)$ are
    orthogonal, so since both operators $d\uPhi_n$ and $D_{n+1}$ are $K$-equivariant
    and linear it suffices to prove equation \eqref{eq:int_equality_D_n_general} for
    elements in $\Gamma(\ch_n)_{\mv}$. Furthermore, we can decompose $\alpha$ into
    $(p,q)$-types and use linearity to expand both sides of equation
    \eqref{eq:int_equality_D_n_general}. Then the claimed equation is true for
    $\alpha$ if and only if it is valid for $\alpha_{p,q}$ for all $p+q =
    n$. Therefore, we fix $0 \le p, q \le n$ with $p+q = n$ and let $\alpha \in
    \Gamma(\ch_{p,q})_{\mv}$ and $\beta \in \Gamma(\ch_n)_{\mv}$. Comparing the types
    of forms on the left and right hand side of \eqref{eq:int_equality_D_n_general}
    it follows that in the decomposition of $\beta$ into $(p,q)$-types only
    $\beta_{p+1,q-1}$, $\beta_{p,q}$ and $\beta_{p-1,q+1}$ can contribute
    nontrivially. Using linearity of $d\uPhi_n$, $D_{n+1}$ and the integral it
    suffices to show the validity of the equation for these three cases. Moreover,
    for $\beta_{p,q}$ the claim was already proven in Theorem
    \ref{thm_relation_integral_D_n}.
  
  For $\beta = \beta_{p+1,q-1} \neq 0$ the left hand side of
  \eqref{eq:int_equality_D_n_general} reduces to the integral over the wedge product
  of $\partial\uPhi_n(\alpha)$ and $\partial \uPhi_n(\overline{\beta}) $, whereas via
  the explicit formula for $D_{n+1}$ from Theorem \ref{thm_BGG_operator}, the right
  hand side is a multiple of $\langle\!\langle \alpha,
  \opartial_H\partial_H^*\beta\rangle\!\rangle = \langle\!\langle
  \opartial_H^*\alpha, \partial_H^*\beta \rangle\!\rangle$. Now if
  $\opartial\uPhi_n(\alpha) = 0$, then by Corollary \ref{cor_Poisson_derivative}, we
  have $\alpha \in \ker(\ocd^*) = \ker(\opartial_H^*)$, showing that both sides of
  \eqref{eq:int_equality_D_n_general} are trivial. By symmetry the analogous
  statement is true if $\partial \uPhi_n(\overline{\beta}) = 0$, so we can assume
  that $\overline{\partial}\uPhi_n(\alpha)$ and $\partial\uPhi_n(\overline{\beta})$ are both
  nontrivial. Now if we are in the first case of Lemma \ref{lem_cases_pairs_degree},
  i.e. if $\alpha$ and $\beta$ are both contained in $\ker(\cd^*) \cap \ker(\ocd^*)$,
  then the exterior derivatives are orthogonal and thus the left hand side of
  \eqref{eq:int_equality_D_n_general} is trivial, whereas the right hand side is
  trivial since $\partial_H^*\alpha = 0$. If we are in the second case of Lemma
  \ref{lem_cases_pairs_degree} can write $\beta = \cd \gamma$ with $\gamma \in
  \Gamma(\ch_{p, q-1})_{\mv}$. Applying the relation $\partial \circ \uPhi_{n-1} =
  \uPhi_n \circ \cd$ and anticommutativity of $\partial$ and $\opartial$ we obtain
  $\partial\uPhi_n(\overline{\beta}) = - \opartial\uPhi_n(\cd \overline{\gamma})$,
  whereas $D_{n+1}(\overline{\beta}) = - D_{n+1}(\cd \overline{\gamma})$ due to
  $D_{n+1} \circ D_n = 0$. Since $\cd\gamma \in \Gamma(\ch_{p,q})_{\mv}$ equation
  \eqref{eq:int_equality_D_n_general} follows again from Theorem
  \ref{thm_relation_integral_D_n}. In the case $\beta = \beta_{p-1,q+1}$ we can argue
  analogously.

(iii) By part (ii), $d\uPhi_n$ has values in
    $L^2\Omega^{n+1}_{\Delta}(G/K)$ and by part (i), it descends to
    $\Gamma(\ch_n)/\im(D_n)$. Thus we only have to show that the image of
    $d\uPhi_n$ is dense. For $0 \le s, t \le n+1$ with $s+t = n+1$ let $\pi_{s,t}
    \colon \Omega^{n+1}(G/K) \to \Omega^{s,t}(G/K)$ be the $G$-equivariant projection
    onto the $(s,t)$-type. Since the Laplace operator is $G$-equivariant and
    preserves types it follows that the image of a harmonic form on $G/K$ under
    $\pi_{s,t}$ is again harmonic. Similarly, the $L^2$-inner product is
    $G$-invariant and different $(s,t)$-types are orthogonal, implying that
    $\pi_{s,t}$ maps between differential forms with finite $L^2$-norm. Therefore, we
    obtain an induced map $\pi_{s,t} \colon L^2\Omega^{n+1}_{\Delta}(G/K) \to
    L^2\Omega^{s,t}_{\Delta}(G/K)$, which is again $G$-equivariant. Combining this
    with $d\uPhi_n$ yiels a continous $G$-equivariant map $\Gamma(\ch_n)/\im(D_n)\to
    L^2\Omega_{\Delta}^{s,t}(G/K)$, which is nontrivial. Now the target space is an
    irreducible $G$-module due to \cite{pedon_complex}*{Theorem 3.4} and thus the
    image has to be dense.
\end{proof}

Formula \eqref{eq:int_equality_D_n_general} explicitly describes the
  pullback of the $L^2$-inner product on $L^2\Omega^{n+1}_{\Delta}(G/K)$ to
  $\Gamma(\ch_n)$. Hence the direct sum of discrete series representations for $G$
  can be identified with the completion of $\Gamma(\ch_n)/\im(D_{n-1})$ with respect
  to this inner product. Now one can rewrite the right hand side of
  \eqref{eq:int_equality_D_n_general} in terms of the operator
  $\sqrt{D_{n+1}^*D_{n+1}}$ to arrive at a description parallel to the one obtained
  by J.\ Lott in \cite{lott} for the case of $\SO_0(n+1,1)$. This was used to identify
  the sum of discrete series representations of the latter group with a space of
  differential forms lying in the Sobolev space $H^{-1/2}$.

\end{document}